\documentclass[11pt,final]{article}

\usepackage{amsfonts,latexsym} 

\usepackage[dvips]{graphics,epsfig,color}
\usepackage{graphicx}
\usepackage{amssymb}
\usepackage{color}
\usepackage{aeguill}

\newcounter{cst}
\def \ctel#1{C_{\refstepcounter{cst}\label{#1}\thecst}}
\def \cter#1{C_{\ref{#1}}}

\newcounter{cexp}
\def \terml#1{T_{\refstepcounter{cexp}\label{#1}\thecexp}}
\def \termr#1{T_{\ref{#1}}}

\newtheorem{defi}{Definition}[section]
\newtheorem{lemma}{Lemma}[section]
\newtheorem{remark}{Remark}[section]
\newtheorem{theo}{Theorem}[section]

\newenvironment{proof}[1]{\noindent {\sc Proof#1} }{$\square$ }
\renewenvironment{abstract}{\begin{center} {\bf Abstract} \\ \end{center}
}{\smallskip}

\def\dsp{\displaystyle}

\def\be{\begin{equation}}
\def\ee{\end{equation}}

\def\ba{\begin{array}{lllll}}
\def\ea{\end{array}}

\def\beqsys {\be\ba \left \{ \begin{array}{l}}
\def\eeqsys {\end{array} \right . \ea\ee }

\def\beqsysno {\be\ba \left \{ \begin{array}{l}}
\def\eeqsysno {\end{array} \right . \ea\ee}

\def\bu{\bar u}
\def\centers{{\cal P}}
\def\card{{{\rm card}}}

\def\di{\bigtriangleup}
\def\disc{{\cal D}}
\def\diam{\hbox{\rm diam}}
\def\div{{\rm div}}
\def\dr{\partial}
\def\dx{\,{\rm d}x}
\def\dfrontiere{\,{\rm d}\gamma}

\def\edge{\sigma}
\def\edges{{\cal E}}
\def\edgesint{{\cal E}_{{\rm int}}}
\def\edgesext{{\cal E}_{{\rm ext}}}

\def\edgescv{{\cal E}_K}

\def\eps{\varepsilon}

\def\grad{\nabla}
\def\lap{\Delta}

\def\mcv{\meas(K)}
\def\medge{\meas(\edge)}
\def\meas{{\rm m}}
\def\mesh{{\cal M}}

\def\n{\mathbf{n}}

\def\NN{{\cal N}}
\def\O{\Omega}
\def\phi{\varphi}

\def\R{\mathbb{R}}
\def\refe#1{(\ref{#1})}
\def\regul{\hbox{\rm regul}}

\def\size{\hbox{\rm size}}

\def\x{{\bf x}}

\usepackage[color,notref,notcite]{showkeys}
\definecolor{labelkey}{rgb}{0.6,0,1}

\begin{document}

\begin{center} {\huge A mixed finite volume scheme for 
anisotropic diffusion problems on any grid}

\large

\vspace*{0.7cm}

J\'er\^ome Droniou \footnote{D\'epartement de Math\'ematiques, UMR CNRS 5149, CC 051, Universit\'e Montpellier II,
Place Eug\`ene Bataillon, 34095 Montpellier cedex 5, France. email: \texttt{droniou@math.univ-montp2.fr}}
and Robert Eymard \footnote{Laboratoire d'Analyse et de Mathématiques Appliquées, UMR 8050,
Universit\'e de Marne-la-Vall\'ee, 5 boulevard Descartes, Champs-sur-Marne, 77454 Marne-la-Vall\'ee Cedex 2,
France. email:\texttt{Robert.Eymard@univ-mlv.fr}}

\vspace*{0.7cm}

27/03/2006

\end{center}

\begin{abstract}
We present a new finite volume scheme for anisotropic heterogeneous diffusion problems on unstructured irregular
grids, which simultaneously gives an approximation of the solution  and of its gradient. 
The approximate solution is shown to converge to the continuous
one as the size of the mesh tends to 0, and an error estimate is given.
An easy implementation method is then proposed, and 
the efficiency of the scheme is shown on various types of grids
and for various diffusion matrices.
\end{abstract}


\section{Introduction}

The computation of an approximate solution for equations involving
a second order elliptic operator is needed 
in so many physical and engineering areas, where
the efficiency of some discretization methods, such as finite difference,
finite element or finite volume methods, has been proven. 
The use of finite volume methods is particularly popular
in the oil engineering field, since it allows for coupled
physical phenomena in the same grids, for which the conservation of
various extensive quantities appears to be a main feature. 
However, it is more challenging to define convergent finite volume schemes
for second-order elliptic operators on
refined, distorted or irregular grids,  designed for the purpose of another problem. 

\medskip

For example, in the framework of geological basin simulation, the grids
are initially fitted on the geological layers boundaries, 
which is a first reason for the loss of orthogonality. Then, these grids are
modified during the simulation, following the compaction of these layers (see \cite{faille}), thus leading
to irregular grids, as those proposed by \cite{kershaw}. 
As a consequence, it is no longer possible to compute
the fluxes resulting from a finite volume scheme for a second order
operator, by a simple two-point difference across each interface between two neighboring control volumes.
Such a two-point scheme is consistent only in the case of an
isotropic operator, using a grid such that the lines connecting the centers of the control volumes are orthogonal to
the edges of the mesh. The problem of finding
a consistent expression using only a small number of points,
for the finite volume fluxes in the general case of any grid
and any anisotropic second order operator, has led to many
works (see \cite{aav}, \cite{aav1}, \cite{aav2},
\cite{faille} and references therein; see also \cite{LMS}).
A recent finite volume scheme has been proposed
\cite{egh, egh2}, permitting to obtain a convergence property
in the case of an anisotropic heterogeneous diffusion problem on unstructured grids,
which all the same satisfy the above orthogonality condition.
In the case where such an orthogonality condition is not satisfied, 
a classical method is the mixed finite element method which 
also gives an approximation of the fluxes and of the gradient of the
unknown (see \cite{ACWY}, \cite{AWY}, \cite{CCJ}, \cite{YAC} for example,
among a very large literature). Note that, although the Raviart-Thomas basis
is not directly available on control volumes which are not
simplices or regular polyhedra, such a basis can
be built on more general irregular grids. 
In \cite{kuznetsov}, such a construction is completed using decomposition into simplices
and a local elimination of the unknowns at the internal edges. In \cite{dehz} and \cite{eghgrad},
such basis functions are obtained from the resolution of a Neumann elliptic problem
in each grid block. However, it has been observed
that the use of mixed finite element method
could demand high refined grids on some highly heterogeneous and anisotropic cases
(see \cite{lepot} and the numerical results provided in the present paper).
Note that an improvement of the mixed finite element scheme is the 
expanded mixed finite element scheme \cite{chen}, where different discrete
approximations are proposed for the unknown, its gradient and the product
of the diffusion matrix by the gradient of the unknown. However,
this last scheme seems to present the same 
restrictions on the meshes  as the mixed finite element scheme.

\medskip

We thus propose in this paper an original finite volume method, called the mixed finite volume method,
which can be applicable on any type of grids in any space dimension,
with very few restrictions on the shape of the  control volumes. 
The implementation of this
scheme is proven to be easy, and no geometric complex shape functions
have to be computed. Accurate results are obtained on coarse grids
in the case of highly heterogeneous and anisotropic problems on irregular grids. 
In order to show the mathematical
and numerical properties of this scheme,
we study here the following problem: find an approximation of 
$\bar u$,
weak solution to the following problem:
\be\ba\dsp
-\div ( \Lambda \grad \bar u) = f \hbox{ in } \Omega, \\
\bar u=0 \hbox{ on } \partial \Omega,
\ea\label{ellgen}\ee
under the following assumptions:
\be
\O \mbox{ is an open bounded connected polygonal subset of }\R^d,
\ d\ge 1,
\label{hypomega}\ee
\be
\begin{array}{c}
\Lambda:\Omega \to M_d(\R)\hbox{ is a bounded measurable function such that}\\
\mbox{there exists $\alpha_0>0$ satisfying $\Lambda(x) \xi\cdot\xi \ge \alpha_0 |\xi|^2$
for a.e. $x\in \Omega$ and all $\xi\in\R^d$},
\end{array}
\label{hyplambda}\ee
(where $M_d(\R)$ stands for the space of $d\times d$ real matrices) and
\be
f \in L^2(\O).
\label{hypfg}\ee
Thanks to Lax-Milgram theorem, there exists a unique weak solution to \refe{ellgen}
in the sense that $\bar u\in H^1_0(\Omega)$ and the equation is satisfied in the
sense of distributions on $\Omega$.

\medskip

The principle of the mixed finite volume scheme, described in Section \ref{secpenal}, 
is the following. We simultaneously
look for approximations $u_K$ and ${\bf v}_K$ in each control 
volume $K$ of $\bar u$ and $\grad \bar u$,
and find an approximation $F_\edge$ at each edge $\edge$ of the mesh of
$\int_{\edge} \Lambda(x) \grad \bar u(x)\cdot \n_{\edge}\dfrontiere(x)$,
where $\n_{\edge}$ is a unit vector normal to $\edge$. The values  $F_\edge$
must then satisfy the conservation equation in each control volume, 
and consistency relations are imposed on  $u_K$, ${\bf v}_K$ and  $F_\edge$.
After having investigated in Section \ref{sec-discspace} the properties of a space associated
with the scheme, we show in Section \ref{sec-study} that it leads to a linear system which
has one and only one approximate solution $u$, ${\bf v}$ and $F$,
and we provide the mathematical analysis of its convergence
and give an error estimate. In Section  \ref{numimp}, we propose an easy implementation
procedure for the scheme, and we use it for the study of
some numerical examples. We thus obtain acceptable
results on some grids for which it would be complex to use other methods, 
or to which empirical methods apply but no mathematical results
of convergence nor stability have yet been obtained.

\section{Definition of the mixed finite volume scheme and main results}\label{secpenal}

We first present the notion of admissible discretization
of the domain $\Omega$, which is necessary to give the expression
of the mixed finite volume scheme.
 
\begin{defi}\label{adisc}{\bf [Admissible discretization]}
Let $\O$ be an open bounded polygonal subset of $\R^d$ ($d\ge 1$), and 
$\dr \O =  \overline{\O}\setminus\O$ its boundary.
An admissible finite volume discretization of $\O$ is given
by $\disc=(\mesh,\edges,\centers)$, where:

\begin{itemize}
\item $\mesh$ is a finite family of
non empty open polygonal convex disjoint subsets of
$\O$ (the ``control volumes'') such that
$\overline{\O}= \dsp{\cup_{K \in \mesh} \overline{K}}$.

\item $\edges$ is a finite family of disjoint subsets of
$\overline{\O}$ (the ``edges'' of the mesh), such that,
for all $\edge\in\edges$, there exists an affine hyperplane $E$ of $\R^d$ and $K\in\mesh$
with $\edge \subset \dr K \cap E$ and $\edge$ is a non empty open convex subset of $E$.
We assume that, for  all $K \in \mesh$, there exists  a subset $\edgescv$ of $\edges$
such that $\dr K  = \dsp{\cup_{\edge \in \edgescv}}\overline{\edge} $.
We also assume that, for all $\sigma\in\edges$, either $\sigma\subset \partial\Omega$
or $\overline{\sigma}=\overline{K}\cap\overline{L}$ for some
$(K,L)\in\mesh\times\mesh$. 

\item $\centers$ is a family of points of $\O$ indexed by $\mesh$, denoted by $\centers = (\x_K)_{K \in \mesh}$
and such that, for all $K \in \mesh$, $\x_K \in K$.

\end{itemize}
\end{defi}

Some examples of admissible meshes in the sense of the above definition
are shown in Figures \ref{gridt} and \ref{grids} in Section \ref{numimp}.
 
\begin{remark} Though the elements of $\mathcal E_K$ may not be the real
edges of a control volume $K$ (each $\sigma\in\mathcal E_K$ may be only a part of a
full edge, see figure \ref{grids}), we will in the following call
``edges of $K$'' the elements of $\mathcal E_K$.

Notice that we could also cut each intersection $\overline{K}\cap \overline{L}$
into more than one edge. This would not change our theoretical results
but this would lead, for practical implementation, to artificially enlarge
the size of the linear systems to solve, which would decrease the
efficiency of the scheme.
\end{remark}

\begin{remark}  The whole mathematical
study done in this paper applies whatever the choice of the point
$\x_K$ in each $K\in\mesh$. In particular, 
we do not impose any orthogonality condition 
connecting the edges and the points $\x_K$.
However, the magnitude of the numerical error (and, for 
some regular or structured types of mesh, its order) does depend on this
choice.

We could also extend our definition to non-planar edges, under some curvature condition.
In this case, it remains possible to use
the mixed finite volume scheme studied in this paper and to prove its convergence.
\end{remark}

The following notations are used. The measure of a control volume $K$ is
denoted by $\mcv$; the $(d-1)$-dimensional measure of an edge $\edge$
is $\meas(\edge)$. In the case where $\sigma\in\mathcal E$
is such that $\overline{\sigma}=\overline{K}\cap\overline{L}$
for $(K,L)\in\mesh\times\mesh$, we denote $\edge = K|L$. 
For all $\edge\in\edges$, $\x_\edge$ is the barycenter
of $\edge$. If $\sigma\in\mathcal E_K$ then $\n_{K,\sigma}$
is the unit normal to $\sigma$ outward to $K$.
The set of interior (resp. boundary) edges is denoted by $\edgesint$
(resp. $\edgesext$), that is $\edgesint = \{\edge \in \edges;$ $\edge \not \subset
\partial \O \}$ (resp.  $\edgesext = \{\edge \in \edges;$ $\edge \subset
\partial \O \}$). For all $K\in\mesh$, we denote by $\NN_K$ the subset of $\mesh$ of the 
neighboring control volumes (that is, the $L$ such that $\overline{K}\cap\overline{L}$
is an edge of the discretization).

To study the convergence of the scheme, we will need the following two
quantities: the size of the discretization
\[
\size(\disc)= \sup\{\hbox{\rm diam}(K)\,;\; K\in \mesh\}
\]
and the regularity of the discretization
\be
\regul(\disc)= \sup\left\{\max\left(\frac{\diam(K)^d}{\rho_K^d},
\mbox{Card}(\mathcal E_K)\right)\,;\; K\in\mesh\right\}
\label{regul}\ee
where, for $K\in \mesh$, $\rho_K$ is the supremum of the radius of the balls
contained in $K$. Notice that, for all $K\in\mesh$,
\begin{equation}
\diam(K)^d \le \regul(\disc)\rho_K^d\le \frac{\regul(\disc)}{\omega_d}\mcv
\label{maj-je}\end{equation}
where $\omega_d$ is the volume of the unit ball in $\R^d$. Note also that
$\regul(\disc)$ does not increase in a local refinement procedure, which
will allow the scheme to handle such procedures.

\medskip

We now define the mixed finite volume scheme. Let $\disc$ be an admissible discretization of $\O$
in the sense of Definition \ref{adisc}. Denote by $H_\disc$ the set of real functions on $\Omega$
which are constant on each control volume $K\in\mathcal M$ (if $h\in H_\disc$, we let $h_K$
be its value on $K$).

As said in the introduction, the idea is to consider three sets of unknowns, namely 
$u\in H_\disc$ which approximates $\bar u$, $\mathbf{v}\in H_\disc^d$
which approximates $\nabla \bar u$ and a family of real numbers
$F=(F_{K,\sigma})_{K\in\mathcal M\,,\sigma\in\mathcal E_K}$ (we denote by
$\mathcal F_\disc$ the set of such families) which approximates $(\int_\sigma \Lambda(x)\nabla \bar u(x)\cdot
\n_{K,\sigma}\,d\gamma(x))_{K\in\mathcal M\,,\sigma\in\mathcal E_K}$. 

Taking $\nu=(\nu_{K})_{K\in\mathcal M}$ a family of nonnegative numbers,
we define $L_\nu(\disc)$ as the space of $(u,\mathbf{v},F)\in H_\disc\times H_\disc^d\times\mathcal F_\disc$
such that
\be\ba
{\bf v}_K\cdot (\x_\edge-\x_K) + {\bf v}_L\cdot (\x_{L} - \x_\edge) + \nu_K\mcv F_{K,\edge}-
\nu_L\meas(L) F_{L,\edge}= u_L - u_K, \\[0.1cm]
\hspace*{6cm} \forall K\in\mesh,\ \forall L\in\NN_K,\mbox{ with } \sigma=K|L,\\[0.2cm]
{\bf v}_K\cdot (\x_\edge-\x_K) + \nu_K \mcv F_{K,\sigma} = - u_K,\quad
\forall K\in\mesh,\ \forall \edge\in\mathcal E_K\cap \mathcal E_{\rm ext}
\ea\label{condconspen2}\ee
and we define the mixed finite volume scheme as: find $(u,\mathbf{v},F)\in L_\nu(\disc)$ such that
\be\ba
F_{K,\sigma}+F_{L,\sigma}=0,\quad \forall \sigma=K|L\in \mathcal E_{\rm int},
\ea\label{condconspen1}\ee
\be\ba
\dsp \meas(K) \Lambda_K {\bf v}_K = \sum_{\edge\in\edgescv} F_{K,\edge} (\x_\edge -
\x_K),\quad \forall K\in\mesh
\ea\label{condconspen3}\ee
(where $\Lambda_K=\frac{1}{\meas(K)}\int_K\Lambda(x)\dx$) and
\be
- \sum_{\edge\in\edgescv} F_{K,\edge} = \int_K f(x)\dx, \quad \forall K\in\mesh.
\label{mul3pen}\ee

The origin of each of these equations is quite easy to understand.
Since $u$ and $\mathbf{v}$ stand for approximate values of $\bar u$ and
$\nabla\bar u$, equation \refe{condconspen2} simply states, if we assume $\nu_K=0$,
that $\mathbf{v}$ is a discrete gradient of $u$:
it is the discrete counterpart of $u(\x_L)-u(\x_K)=u(\x_L)-u(\x_\sigma)
+u(\x_\sigma)-u(\x_K)\approx \nabla u(\x_L)\cdot(\x_L-\x_\sigma)
+\nabla u(\x_K)\cdot(\x_\sigma-\x_K)$. This equation is slightly penalized with the fluxes to ensure
existence and estimates on the said fluxes (to study the convergence of the scheme,
we will assume $\nu_K>0$; see the theorems below).
Notice that the boundary condition $\bar u=0$ is contained in the second line of \refe{condconspen2}.
As $F_{K,\sigma}$ stands for an approximate value 
of $\int_\sigma \Lambda \nabla(x)\bar u(x)\cdot\n_{K,\sigma}\,d\gamma(x)$,
it is natural to ask for the conservation property \refe{condconspen1},
and the balance \refe{mul3pen}
simply comes from an integration of \refe{ellgen} on a control volume.
Last, the link \refe{condconspen3} between $\Lambda\mathbf{v}$ and its fluxes
is justified by Lemma \ref{form-mag} in the appendix, which shows that one
can reconstruct a vector from its fluxes through the edges of a control volume.

\medskip

Our main results on the mixed finite volume scheme are the following. The first one states that there
exists a unique solution to the scheme. The second one gives the convergence
of this solution to the solution of the continuous problem, as the size
of the mesh tends to $0$, and the third one provides an error estimate
in the case of smooth data.

\begin{theo}\label{existpen}
Let us assume Assumptions \refe{hypomega}-\refe{hypfg}. Let $\disc$ be an admissible
discretization of $\O$ in the sense of Definition \ref{adisc}.
Let $(\nu_K)_{K\in\mesh}$ be a family of positive real numbers.
Then there exists one and only one $(u,{\bf v},F)$ solution of 
(\refe{condconspen2},\refe{condconspen1},\refe{condconspen3},\refe{mul3pen}).
\end{theo}

\begin{theo}\label{convpen}
Let us assume Assumptions \refe{hypomega}-\refe{hypfg}. Let $(\disc_m)_{m\ge 1}$ be
admissible discretizations of $\O$ in the sense of Definition \ref{adisc},
such that $\size(\disc_m)\to 0$ as $m\to\infty$ and $(\regul(\disc_m))_{m\ge 1}$
is bounded. Let $\nu_0>0$ and $\beta\in (2-2d,4-2d)$ be fixed.
For all $m\ge 1$, let $(u_m,\mathbf{v}_m,F_m)$ be the solution to 
(\refe{condconspen2},\refe{condconspen1},\refe{condconspen3},\refe{mul3pen}) for
the discretization $\disc_m$, setting  $\nu_K = \nu_0 \diam(K)^\beta$ for all $K\in\mesh_m$. 
Let $\bar u$ be the weak solution to \refe{ellgen}.

Then, as $m\to\infty$, $\mathbf{v}_m\to \nabla \bar u$ strongly in $L^2(\O)^d$
and $u_m\to\bar u$ weakly in $L^2(\O)$ and strongly in $L^q(\O)$ for all $q<2$.
\end{theo}

\begin{theo}\label{esterrorpen}
Let us assume Assumptions \refe{hypomega}-\refe{hypfg}. Let $\disc$ be an
admissible discretization of $\O$ in the sense of Definition \ref{adisc},
such that $\size(\disc)\le 1$ and $\regul(\disc)\le \theta$ for some $\theta>0$.
We take $\nu_0>0$ and $\beta\in (2-2d,4-2d)$ and, for all $K\in\mathcal M$,
we let $\nu_K = \nu_0 \diam(K)^\beta$.
Let $(u,\mathbf{v},F)$ be the solution to 
(\refe{condconspen2},\refe{condconspen1},\refe{condconspen3},\refe{mul3pen}). 
Let $\bar u$ be the weak solution to
\refe{ellgen}.
We assume that $\Lambda\in C^1(\overline{\Omega};M_d(\R))$ and $\bar u\in C^2(\overline{\Omega})$.

Then there exists $\ctel{pencerr0}$ only depending on $d$, $\O$,
$\bar u$, $\Lambda$, $\theta$ and $\nu_0$ such that
\be
\Vert \mathbf{v} - \grad \bar u\Vert_{L^2(\O)^d} \le \cter{pencerr0}
\size(\disc)^{\frac{1}{2}\min(\beta +2d-2,4-2d - \beta )}
\label{peneqerr0}\ee
and
\be
\Vert u - \bar u\Vert_{L^2(\O)} \le \cter{pencerr0}
\size(\disc)^{\frac{1}{2}\min(\beta +2d -2,4-2d - \beta )}
\label{peneqerr1}\ee
(note that the maximum value of $\frac{1}{2}\min(\beta +2d-2,4-2d - \beta )$
is $\frac{1}{2}$, obtained for $\beta = 3 - 2d$).
\end{theo} 

\begin{remark} These error estimates are not sharp, and the numerical results
in Section \ref{numimp} show a much better order of convergence.
\end{remark}

\section{The discretization space}\label{sec-discspace}

We investigate here some properties of the space $L_\nu(\disc)$, which will be
useful to study the mixed finite volume
scheme. Recall that $L_\nu(\disc)$ is the space of $(u,\mathbf{v},F)$
which satisfy \refe{condconspen2}.

\begin{lemma} {\rm\bf [Poincar\'e's Inequality]}\label{pdiscpen}
Let us assume Assumption \refe{hypomega}. Let $\disc$ be an admissible
discretization of $\O$ in the sense of Definition \ref{adisc}, such that
$\regul(\disc)\le \theta$ for some $\theta>0$.
Let $(\nu_K)_{K\in\mesh}$ be a family of nonnegative real numbers.
Then there exists $\ctel{poinpen}$ only depending on $d$, $\O$ and $\theta$
such that, for all $(u,{\bf v},F)\in L_\nu(\disc)$,
\be
\Vert u\Vert_{L^2(\O)}\le 
\cter{poinpen}\left(\Vert {\bf v}\Vert_{L^2(\O)^d}+N_2(\disc,\nu,F)\right),
\label{pdispen}\ee
where we have noted $N_2(\disc,\nu,F)=\left(\sum_{K\in\mesh}
\sum_{\edge\in\edgescv} \diam(K)^{2d-2}\nu_K^2 F_{K,\sigma}^2 \mcv\right)^{1/2}$.
\end{lemma}

\begin{proof}{.}

Let $R>0$ and $x_0\in\O$ be such that $\O\subset B(x_0,R)$ (the open ball
of center $x_0$ and radius $R$). We extend $u$ by the value $0$
in $B(x_0,R)\setminus\O$, and we consider $w\in H^1_0(B(x_0,R))\cap H^2(B(x_0,R))$ such that
$-\lap w(x) = u(x)$ for a.e. $x\in B(x_0,R)$. We multiply each equation of \refe{condconspen2}
by $\int_{\edge} \grad w(x)\cdot \n_{K,\edge}\dfrontiere(x)$, and we sum on $\edge\in\edges$;
since $\n_{K,\sigma}=-\n_{L,\sigma}$ whenever $\sigma=K|L$, we find
\begin{eqnarray*}
&&\sum_{\sigma\in\mathcal E_{\rm int}\,,\sigma=K|L} \mathbf{v}_K\cdot(\x_\sigma-\x_K)\int_{\edge} \grad w(x)\cdot \n_{K,\edge}\dfrontiere(x)
+\mathbf{v}_L\cdot(\x_\sigma-\x_L)\int_{\edge} \grad w(x)\cdot \n_{L,\edge}\dfrontiere(x)\\
&&+\sum_{\sigma\in\mathcal E_{\rm ext}\,,\sigma\in\mathcal E_K}\mathbf{v}_K\cdot(\x_\sigma-\x_K)\int_{\edge} \grad w(x)\cdot \n_{K,\edge}\dfrontiere(x)\\
&&+\sum_{\sigma\in\mathcal E_{\rm int}\,,\sigma=K|L} \nu_K \meas(K)F_{K,\sigma}\int_{\edge} \grad w(x)\cdot \n_{K,\edge}\dfrontiere(x)
+\nu_L\meas(L)F_{L,\sigma}\int_{\edge} \grad w(x)\cdot \n_{L,\edge}\dfrontiere(x)\\
&&+\sum_{\sigma\in\mathcal E_{\rm ext}\,,\sigma\in\mathcal E_K}
\nu_K\meas(K)F_{K,\sigma}\int_{\edge} \grad w(x)\cdot \n_{K,\edge}\dfrontiere(x)\\
&=&
-\sum_{\sigma\in\mathcal E_{\rm int}\,,\sigma=K|L} u_K\int_{\edge} \grad w(x)\cdot \n_{K,\edge}\dfrontiere(x)
+u_L\int_{\edge} \grad w(x)\cdot \n_{L,\edge}\dfrontiere(x)\\
&&-\sum_{\sigma\in\mathcal E_{\rm ext}\,,\sigma\in\mathcal E_K}u_K\int_{\edge} \grad w(x)\cdot \n_{K,\edge}\dfrontiere(x).
\end{eqnarray*}

Gathering by control volumes, we find
\begin{eqnarray}
\sum_{K\in\mesh}{\bf v}_K\cdot 
\sum_{\edge\in\edgescv}(\x_\edge-\x_K)
\int_{\edge} \grad w(x)\cdot \n_{K,\edge} \dfrontiere(x)\nonumber&&\\
+\sum_{K\in\mesh}
\sum_{\edge\in\edgescv} \nu_K \mcv F_{K,\sigma}
\int_{\edge} \grad w(x)\cdot \n_{K,\edge} \dfrontiere(x)
&=&-\sum_{K\in\mesh} u_K \sum_{\sigma\in\mathcal E_K} \int_{\edge} 
\grad w(x)\cdot \n_{K,\edge}\dfrontiere(x)\nonumber\\
&=&-\sum_{K\in\mesh} u_K \int_{K} \Delta w(x)\,dx\nonumber\\
&=&\sum_{K\in\mesh} \mcv u_K^2=||u||_{L^2(\O)}^2.\label{rajout}
\end{eqnarray}
Let $\terml{pdp1}$ and $\terml{pdp2}$ be the two terms in the left-hand
side of this equation.

\medskip

Define $\terml{pdp1b} = \int_\O {\bf v}(x) \cdot \grad w(x) \dx$; we have
\begin{equation}
|\termr{pdp1b}| \le \Vert {\bf v}\Vert_{L^2(\O)^d}\Vert w \Vert_{H^1(\O)}
\label{rajout0}\end{equation}
and we want to compare $\termr{pdp1}$ with $\termr{pdp1b}$.
In order to do so, we apply Lemma \ref{form-mag} in the appendix to the vector
${\bf G}_K = \frac 1 \mcv \int_K \grad w(x) \dx$, which gives
\[
\int_K \nabla w(x)\dx=\mcv \mathbf{G}_K=\sum_{\sigma\in\mathcal E_K}
\meas(\sigma) \mathbf{G}_K\cdot\n_{K,\sigma}(\x_\edge-\x_K)
\]
and therefore
\[
\termr{pdp1b} = \sum_{K\in\mesh}{\bf v}_K\cdot 
\sum_{\edge\in\edgescv}\medge
{\bf G}_K \cdot \n_{K,\edge}(\x_\edge-\x_K).
\]
Hence, setting ${\bf G}_{\edge} = \frac 1 \medge \int_\edge 
\grad w(x) \dfrontiere(x)$, we get
\[
|\termr{pdp1} - \termr{pdp1b}| \le \sum_{K\in\mesh}|{\bf v}_K|
\sum_{\edge\in\edgescv}\medge\
|{\bf G}_K  - {\bf G}_{\edge}|\ \diam(K).
\]
Thanks to the Cauchy-Schwarz inequality, we find
\[
(\termr{pdp1} - \termr{pdp1b})^2 \le \left( \sum_{K\in\mesh} |{\bf v}_K|^2
\sum_{\edge\in\edgescv}\medge\diam(K)\right)\left( 
\sum_{K\in\mesh} \sum_{\edge\in\edgescv}\medge\diam(K)
|{\bf G}_K  - {\bf G}_{\edge}|^2\right).
\]
We now apply Lemma \ref{est-intervf} in the appendix, which gives $\ctel{averifier}$ only depending
on $d$ and $\theta$ such that
\begin{equation}
|{\bf G}_K  - {\bf G}_{\edge}|^2 \le \cter{averifier} \frac {\diam(K)} \medge
\Vert w \Vert_{H^2(K)}^2
\label{maj-G}\end{equation}
(notice that $\alpha:=\frac{1}{2} \theta^{-1/d} < \regul(\disc)^{-1/d}
\le \rho_K/\diam(K)$ is valid in Lemma \ref{est-intervf}).
We also have, for $\sigma\in\mathcal E_K$, $\meas(\sigma)\le \omega_{d-1}\diam(K)^{d-1}$,
where $\omega_{d-1}$ is the volume of the unit ball in $\R^{d-1}$.
Therefore, according to \refe{maj-je} and since
$\regul(\disc)\ge \card(\edgescv)$ for all $K\in\mesh$,
\begin{eqnarray}
(\termr{pdp1} - \termr{pdp1b})^2 &\le& \left( \sum_{K\in\mesh} |{\bf v}_K|^2
\sum_{\edge\in\edgescv}\medge\ \diam(K)\right)
\left( \sum_{K\in\mesh}\sum_{\sigma\in\mathcal E_K} \cter{averifier}\diam(K)^2 ||w||_{H^2(K)}^2\right)\nonumber\\
&\le& \left( \omega_{d-1}\regul(\disc)\sum_{K\in\mesh} |{\bf v}_K|^2 \diam(K)^d\right)
\left(\cter{averifier} \size(\disc)^2\regul(\disc)\Vert w \Vert_{H^2(\O)}^2 \right)\nonumber\\
&\le& \frac{\omega_{d-1}\regul(\disc)^2}{\omega_d} ||\mathbf{v}||_{L^2(\O)^d}^2
\cter{averifier} \diam(\O)^2\regul(\disc)\Vert w \Vert_{H^2(\O)}^2.
\label{rajout2}
\end{eqnarray}

\medskip

Turning to $\termr{pdp2}$, we have $\termr{pdp2}=\sum_{K\in\mesh}\sum_{\edge\in\edgescv} \nu_K \mcv F_{K,\sigma}
\meas(\sigma)\mathbf{G}_{\sigma}\cdot\n_{K,\sigma}$,
which we compare with $\terml{pdp2b} = \sum_{K\in\mesh}\sum_{\edge\in\edgescv} \nu_K  \mcv F_{K,\sigma} \medge
{\bf G}_K\cdot \n_{K,\edge}$
thanks to \refe{maj-G}:
\begin{eqnarray}
(\termr{pdp2}\lefteqn{-\termr{pdp2b})^2}&&\nonumber\\
&\le& \left(\sum_{K\in\mesh}
\sum_{\edge\in\edgescv} \diam(K)\medge\nu_K^2 F_{K,\sigma}^2 \mcv^2\right)
\left(\sum_{K\in\mesh}\sum_{\edge\in\edgescv}\frac{\medge}{\diam(K)}
|{\bf G}_K  - {\bf G}_{\edge}|^2\right)\nonumber\\
&\le& \left(\omega_{d-1}\omega_d\sum_{K\in\mesh}
\sum_{\edge\in\edgescv} \diam(K)^{2d}\nu_K^2 F_{K,\sigma}^2 \mcv\right)
\regul(\disc)\cter{averifier}||w||_{H^2(\O)}^2\nonumber\\
&\le& \omega_{d-1}\omega_d \diam(\O)^2 N_2(\disc,\nu,F)^2
\regul(\disc)\cter{averifier}||w||_{H^2(\O)}^2.
\label{rajout3}
\end{eqnarray}
On the other hand, we can write
\begin{eqnarray}
\termr{pdp2b}^2 &\le& \left(\sum_{K\in\mesh}
\sum_{\edge\in\edgescv} \medge^2\nu_K^2 F_{K,\sigma}^2 \mcv\right)\left(\sum_{K\in\mesh}
\sum_{\edge\in\edgescv} \mcv|{\bf G}_K|^2\right)\nonumber\\
&\le& \omega_{d-1}^2 N_2(\disc,\nu,F)^2
\left(\regul(\disc)\sum_{K\in\mesh} \mcv |{\bf G}_K|^2\right)\nonumber\\
&\le& \omega_{d-1}^2 N_2(\disc,\nu,F)^2\regul(\disc)\Vert w \Vert_{H^1(\O)}^2.
\label{rajout4}
\end{eqnarray}

\medskip

Thanks to \refe{rajout0}, \refe{rajout2}, \refe{rajout3} and \refe{rajout4}, we
can come back in \refe{rajout} to find
\begin{eqnarray*}
||u||_{L^2(\O)}^2&=&\termr{pdp1} + \termr{pdp2}\\
        &\le& |\termr{pdp1}-\termr{pdp1b}|+|\termr{pdp1b}|+|\termr{pdp2}-\termr{pdp2b}|+|\termr{pdp2b}|\\
&\le& \sqrt{\frac{\omega_{d-1}\cter{averifier}\theta^3}{\omega_d}}\diam(\O)
\Vert \mathbf{v} \Vert_{L^2(\O)^d}||w||_{H^2(\O)}
+||\mathbf{v}||_{L^2(\O)^d}||w||_{H^1(\O)}\\
&&+\sqrt{\omega_{d-1}\omega_d\cter{averifier}\theta}\,\diam(\O)N_2(\disc,\nu,F)\Vert w \Vert_{H^2(\O)}
+\omega_{d-1}\sqrt{\theta}\,N_2(\disc,\nu,F)\Vert w \Vert_{H^1(\O)}.
\end{eqnarray*}
Since there exists $\ctel{reg}$ only depending on $d$ and $B(x_0,R)$ (the ball chosen
at the beginning of the proof) such that $||w||_{H^2(\O)}\le \cter{reg}||u||_{L^2(\O)}$,
this concludes the proof. \end{proof}

\begin{lemma} {\rm\bf [Equicontinuity of the translations]}\label{pdistrpen}
Let us assume Assumption \refe{hypomega}. Let $\disc$ be an admissible
discretization of $\O$ in the sense of Definition \ref{adisc}, such that
$\regul(\disc)\le \theta$ for some $\theta>0$.
Let $(\nu_K)_{K\in\mesh}$ be a family of nonnegative real numbers.
Then there exists $\ctel{poinpen2}$ only depending on $d$, $\O$ and $\theta$
such that, for all $(u,{\bf v},F)\in L_\nu(\disc)$ and all $\xi\in\R^d$,
\be
\Vert u(\cdot+\xi) - u\Vert_{L^1({\R^d})} \le 
\cter{poinpen2}\left(\Vert {\bf v}\Vert_{L^1(\O)^d} + N_1(\disc,\nu,F)\right)|\xi|,
\label{translatpen}\ee
where $N_1(\disc,\nu,F)=\sum_{K\in\mesh}\sum_{\edge\in\edgescv}
\diam(K)^{d-1}\nu_K|F_{K,\edge}| \mcv$ (and $u$ has been extended by
$0$ outside $\O$).
\end{lemma} 

\begin{proof}{.}

For all $\edge\in\edges$, let us define $D_\edge u  = |u_L - u_K|$ if $\edge=K|L$ and
$D_\edge u  = |u_K|$ if $\edge\in\mathcal E_K\cap \mathcal E_{\rm ext}$.  For $(x,\xi)\in \R^d\times\R^d$
and $\sigma\in\edges$, we define $\chi(x,\xi,\edge)$ by $1$ if 
$\edge\cap [x,x+\xi]\neq\emptyset$ and by $0$ otherwise.
We then have, for all $\xi\in\R^d$ and a.e. $x\in\R^d$ (the $x$'s such that
$x$ and $x+\xi$ do not belong to $\cup_{K\in\mesh} \partial K$,
and $[x,x+\xi]$ does not intersect the relative boundary of any edge),
\[
|u(x+\xi) - u(x)| \le \sum_{\edge\in\edges} \chi(x,\xi,\edge) D_\edge u.
\]
Applying \refe{condconspen2}, we get $|u(x+\xi) - u(x)| \le \terml{tp1}(x) +\terml{tp2}(x)$
with
\begin{eqnarray*}
\termr{tp1}(x) &=& \sum_{\edge\in\edgesint,\edge=K|L} \chi(x,\xi,\edge) 
(|{\bf v}_K| |\x_\edge-\x_K| + |{\bf v}_L| |\x_L - \x_\edge|) \\
&&+\sum_{\edge\in\edgesext,\edge\in\edgescv} \chi(x,\xi,\edge) 
|{\bf v}_K| |\x_\edge-\x_K|\\
&\le&\sum_{K\in\mesh}\sum_{\sigma\in\mathcal E_K} \chi(x,\xi,\sigma) \diam(K)|{\bf v}_K|
\end{eqnarray*}
and 
\begin{eqnarray*}
\termr{tp2}(x) &=& \sum_{\edge\in\edgesint,\edge=K|L} 
\chi(x,\xi,\edge) \left(\nu_K \meas(K) |F_{K,\edge}|+ \nu_L\meas(L)|F_{L,\edge}|\right) \\
&&+\sum_{\edge\in\edgesext,\edge\in\edgescv} \chi(x,\xi,\edge) \nu_K  \meas(K) |F_{K,\edge}|\\
&=&\sum_{K\in\mesh}\sum_{\edge\in\edgescv}\chi(x,\xi,\edge)\nu_K \mcv |F_{K,\edge}|.
\end{eqnarray*}

In order that $\chi(x,\xi,\sigma)\not=0$, $x$ must lie in the set $\sigma-[0,1]\xi$
which has measure $\mbox{m}(\sigma)|\mathbf{n}_\sigma\cdot\xi|$
(where $\mathbf{n}_\sigma$ is a unit normal to $\sigma$). Hence,
\[
\int_{\R^d} \chi(x,\xi,\sigma)\,dx\le
\mbox{m}(\sigma)|\mathbf{n}_\sigma\cdot\xi|
\le \omega_{d-1}\diam(K)^{d-1}|\xi|\qquad\mbox{ if $\sigma\in\mathcal E_K$}.
\]
Since ${\rm Card}(\mathcal E_K)\le \regul(\disc)$, this gives
\[
\int_{\R^d} \termr{tp1}(x)\,dx\le \omega_{d-1}\regul(\disc) |\xi|\sum_{K\in\mesh} \diam(K)^d |\mathbf{v}_K|
\]
and
\[
\int_{\R^d} \termr{tp2}(x)\dx \le\omega_{d-1}|\xi| \sum_{K\in\mesh}\sum_{\edge\in\edgescv}
\diam(K)^{d-1}\nu_K|F_{K,\edge}| \mcv,
\]
which concludes the proof thanks to \refe{maj-je}. \end{proof}

\begin{remark}
We could prove the property 
\[
\Vert u(\cdot+\xi) - u\Vert_{L^2(\R^d)}^2 \le 
C (\Vert {\bf v}\Vert_{L^2(\O)^d}^2+N_2(\disc,\nu,F)^2)\ |\xi|\ (|\xi| + \size(\disc)),
\]
by introducing the maximum value of $\diam(K)/\rho_L$, 
for all $(K,L)\in\mesh\times\mesh$, in the definition of $\regul(\disc)$. This
would allow, in Theorem \ref{convpen}, to prove the strong convergence of $u_m$
in $L^2(\O)$. Nevertheless, we chose not to do so since
the quantity $\diam(K)/\rho_L$ cannot remain bounded in a local mesh refinement
procedure. 

Note that we shall all the same prove, in a particular case, 
a strong convergence property in $L^2(\O)$ for $u$, as a consequence of the error
estimate. 
\end{remark}

\begin{lemma} {\rm\bf [Compactness property]}\label{cpctpen}
Let us assume Assumption \refe{hypomega}. Let $(\disc_m)_{m\ge 1}$ be admissible
discretizations of $\O$ in the sense of Definition \ref{adisc}, such that
$\size(\disc_m)\to 0$ as $m\to\infty$ and $(\regul(\disc_m))_{m\ge 1}$ is bounded.
Let $(u_m,{\bf v}_m,F_m,\nu_m)_{m\ge 1}$ be such that
$(u_m,{\bf v}_m,F_m)\in L_{\nu_m}(\disc_m)$, $(\mathbf{v}_m)_{m\ge 1}$
is bounded in $L^2(\O)^d$ and $N_2(\disc_m,\nu_m,F_m)\to 0$ as $m\to\infty$
($N_2$ has been defined in Lemma \ref{pdiscpen}).

Then there exists a subsequence of $(\disc_m)_{m\ge 1}$ (still denoted by $(\disc_m)_{m\ge 1}$)
and $\bar u\in H^1_0(\O)$ such that the corresponding sequence
$(u_m)_{m\ge 1}$ converges to $\bar u$ weakly in $L^2(\O)$ and strongly in $L^q(\O)$ for
all $q<2$, and such that $({\bf v}_m)_{m\ge 1}$ converges to $\grad\bar u$ weakly in $L^2(\O)^d$.
\end{lemma} 

\begin{proof}{.}

Notice first that, for all discretization $\disc$, for all $\nu=(\nu_K)_{K\in\mesh}$
nonnegative numbers and for all $F=(F_{K,\sigma})_{K\in\mesh\,,\;\sigma\in\mathcal E_K}$,
\begin{eqnarray}
N_1(\disc,\nu,F)&=&\sum_{K\in\mesh}\sum_{\sigma\in\mathcal E_K}\diam(K)^{d-1}\nu_K |F_{K,\sigma}|\mcv\nonumber\\
&\le& \left(\sum_{K\in\mesh}\sum_{\sigma\in\mathcal E_K}\diam(K)^{2d-2}\nu_K^2 F_{K,\sigma}^2
\mcv\right)^{1/2} \left(\sum_{K\in\mesh}\sum_{\sigma\in\mathcal E_K}\mcv\right)^{1/2}\nonumber\\
&\le& N_2(\disc,\nu,F) \regul(\disc)^{1/2}\meas(\Omega)^{1/2}.
\label{debutproof}\end{eqnarray}
Hence, if $N_2(\disc,\nu,F)$ and $\regul(\disc)$ are bounded, so is $N_1(\disc,\nu,F)$.
Owing to this, the hypotheses, Lemmas \ref{pdiscpen}, \ref{pdistrpen} 
and Kolmogorov's compactness theorem allow to extract
a subsequence such that $\mathbf{v_m}\to \bar \mathbf{v}$ weakly in $L^2(\O)^d$
and $u_m\to \bar u$ weakly in $L^2(\O)$ and strongly in $L^1(\Omega)$ (which
implies the strong convergence in $L^q(\O)$ for all $q < 2$).
We now extend $u_m$, $\bar u$, $\mathbf{v}_m$ and $\bar\mathbf{v}$ by $0$ outside
$\O$ and we prove that $\bar{\mathbf{v}}=\nabla \bar u$ in the distributional
sense on $\R^d$. This will conclude that $\bar u\in H^1(\R^d)$ and, since
$\bar u=0$ outside $\O$, that $\bar u\in H^1_0(\O)$.

\medskip

Let ${\bf e}\in\R^d$ and $\varphi\in C^\infty_c(\R^d)$. For simplicity, we drop the index $m$
for $\disc_m$, ${\bf v}_m$ and $u_m$. We multiply each equation of
\refe{condconspen2} by $\int_{\sigma} \varphi(x) \dfrontiere(x) {\bf e}\cdot \n_{K,\sigma}$,
we sum all these equations and we gather by control volumes, getting
$\terml{t1} +\terml{t11}= \terml{t2}$ with
\[
\termr{t1} = \sum_{K\in\mesh} {\bf v}_K\cdot \sum_{\edge\in\edgescv} 
\int_{\edge} \varphi(x) \dfrontiere(x)\; {\bf e}\cdot \n_{K,\edge}(\x_\edge-\x_K),
\]
\[
\termr{t11} = \sum_{K\in\mesh} \sum_{\sigma\in\mathcal E_K}\nu_K \meas(K)F_{K,\sigma}
\int_{\sigma} \varphi(x) \dfrontiere(x) {\bf e}\cdot \n_{K,\sigma}
\]
and 
\[
\termr{t2} =- \sum_{K\in\mesh} u_K \sum_{\edge\in\edgescv} 
\int_{\edge} \varphi(x) \dfrontiere(x)\ {\bf e}\cdot \n_{K,\edge} = 
- \int_\O u(x) \div  (\varphi(x) {\bf e})\dx.
\]

\medskip

We want to compare $\termr{t1}$ with $\terml{t1b}$ defined by
\[
\termr{t1b} = \sum_{K\in\mesh} {\bf v}_K\cdot \sum_{\edge\in\edgescv} 
\frac 1 {\mcv} \int_{K} \varphi(x) \dx\; \medge \  {\bf e}\cdot
\n_{K,\edge}(\x_\edge-\x_K).
\]
Since there exists $\ctel{varphi}$ only depending on $\varphi$ such that, for all
$\sigma\in\mathcal E_K$,
\[
\left\vert \frac 1 {\medge}\int_{\edge} \varphi(x) \dfrontiere(x) - 
\frac 1 {\mcv} \int_{K} \varphi(x) \dx\right\vert \le \cter{varphi}\size(\disc),
\]
we get that
\[
\left\vert\termr{t1} - \termr{t1b}\right\vert \le \cter{varphi}|\mathbf{e}|\size(\disc)
\sum_{K\in\mesh} \vert{\bf v}_K\vert \sum_{\edge\in\edgescv}  
\medge \vert\x_\edge-\x_K\vert.
\]
But $\medge\vert\x_\edge-\x_K\vert \le \omega_{d-1}\diam(K)^d\le
\frac{\omega_{d-1}{\rm regul}(\disc)}{\omega_d}\mcv$ and, since
$\card(\edgescv)\le \regul(\disc)$, we obtain
\[
\left\vert\termr{t1} - \termr{t1b}\right\vert \le
\cter{varphi}|\mathbf{e}|\size(\disc)\frac{\omega_{d-1}\regul(\disc)^2}{\omega_d}
\Vert \mathbf{v}\Vert_{L^1(\O)}
\]
and thus $\lim_{\mbox{\small size}(\disc)\to 0} |\termr{t1} - \termr{t1b}| =  0$.
Moreover, thanks to Lemma \ref{form-mag}, we get $\termr{t1b} = \int_\O \varphi(x)  
{\bf v}(x) \cdot {\bf e} \dx$
and so $\lim_{\mbox{\small size}(\disc)\to 0} \termr{t1b} =  \int_\O \varphi(x)  
\bar {\bf v}(x) \cdot {\bf e} \dx=\int_{\R^d} \varphi(x) \bar {\bf v}(x) \cdot {\bf e} \dx$
($\bar {\bf v}$ has been extended by $0$ outside $\O$).
This proves that
\begin{equation}
\lim_{\mbox{\small size}(\disc)\to 0}\termr{t1}=
\int_{\R^d} \varphi(x)  \bar {\bf v}(x) \cdot {\bf e} \dx.
\label{rajout112}\end{equation}

\medskip

Since $\varphi$ is bounded, by \refe{debutproof} we find $\ctel{varphie}$ only depending
on $\varphi$ and $\mathbf{e}$ such that
\begin{eqnarray*}
|\termr{t11}|&\le& \cter{varphie} \sum_{K\in\mesh} \sum_{\edge\in\edgescv} \medge \nu_K |F_{K,\edge}| \mcv\\
&\le& \cter{varphie}\omega_{d-1}N_1(\disc,\nu,F)\\
&\le& \cter{varphie}\omega_{d-1}\regul(\disc)^{1/2}\meas(\Omega)^{1/2}N_2(\disc,\nu,F)
\end{eqnarray*}
and therefore, by the assumptions,
\begin{equation}
\lim_{\mbox{\small size}(\disc)\to 0}\termr{t11}=0.
\label{rajout122}\end{equation}

We clearly have
\[
\lim_{\mbox{\small size}(\disc)\to 0} \termr{t2} = - \int_\O \bar u(x) \div  (\varphi(x) {\bf e})\dx
= - \int_{\R^d} \bar u(x) \div  (\varphi(x) {\bf e})\dx
\]
(recall that $\bar u$ has been extended by $0$ outside $\O$). Gathering
this limit with \refe{rajout112} and \refe{rajout122} in $\termr{t1}+\termr{t11}=\termr{t2}$,
we obtain
\[
\int_{\R^d} \varphi(x)\bar\mathbf{v}(x)\cdot \mathbf{e}\,dx=
-\int_{\R^d} \bar u(x) \div  (\varphi(x) {\bf e})\dx,
\]
which concludes the proof that $\bar \mathbf{v}=\nabla \bar u$ in the
distributional sense on $\R^d$. \end{proof}

\section{Study of the mixed finite volume scheme}\label{sec-study}

We first prove an \emph{a priori} estimate on the solution to 
the scheme.

\begin{lemma}\label{exestpen}
Let us assume Assumptions \refe{hypomega}-\refe{hypfg}. Let $\disc$ be an admissible
discretization of $\O$ in the sense of Definition \ref{adisc}.
Let $(\nu_K)_{K\in\mesh}$ be a family of positive real numbers
and $(u,{\bf v},F)\in L_\nu(\disc)$ be a solution of 
(\refe{condconspen2},\refe{condconspen1},\refe{condconspen3},\refe{mul3pen}).
Then, for all $\nu_0 >0$, for all $\beta_0\ge \beta \ge 2-2d$ such that
$\nu_K \le \nu_0 \diam(K)^{\beta}$ ($\forall K\in\mesh$) and for all $\theta\ge \regul(\disc)$, this solution satisfies
\begin{equation}
\Vert{\bf v}\Vert_{L^2(\O)^d}^2+\sum_{K\in\mesh}\sum_{\sigma\in\mathcal E_K}
\nu_K F_{K,\sigma}^2\mcv\le \ctel{poinpen-je}||f||_{L^2(\O)}^2
\label{estimpen}\end{equation}
where $\cter{poinpen-je}$ only depends on $d$, $\O$, $\alpha_0$, $\theta$, $\nu_0$ and $\beta_0$.
\end{lemma} 

\begin{remark} This estimate shows that if $f=0$ then $F=0$ and $\mathbf{v}=0$, and thus $u=0$ by Lemma \ref{pdiscpen}.
Since (\refe{condconspen2},\refe{condconspen1},\refe{condconspen3},\refe{mul3pen}) is square and
linear in $(u,\mathbf{v},F)$, the existence and uniqueness of the solution to the mixed
finite volume scheme (i.e. Theorem \ref{existpen}) is an immediate consequence of this lemma.
\label{proofth}\end{remark}

\begin{proof}{.}

Multiply \refe{mul3pen} by $u_K$, sum on the control volumes and gather by edges using \refe{condconspen1}:
\[
\sum_{\sigma\in\mathcal E_{\rm int}\,,\sigma=K|L}F_{K,\sigma}(u_L-u_K)+
\sum_{\sigma\in\mathcal E_{\rm ext}\,,\sigma\in\mathcal E_K}-F_{K,\sigma}u_K
=\int_\Omega f(x)u(x)\,dx.
\]
Using \refe{condconspen2} and \refe{condconspen1}, and gathering by control volumes,
this gives
\begin{eqnarray*}
\lefteqn{\int_\Omega f(x)u(x)\,dx}&&\\
&=&\sum_{\sigma\in\mathcal E_{\rm int}\,,\sigma=K|L}F_{K,\sigma}\mathbf{v}_K\cdot(\x_\sigma-\x_K)+
F_{L,\sigma}\mathbf{v}_L\cdot(\x_\sigma-\x_L)+
\sum_{\sigma\in\mathcal E_{\rm ext}\,,\sigma\in\mathcal E_K}F_{K,\sigma}\mathbf{v}_K\cdot(\x_\sigma-\x_K)\\
&&+\sum_{\sigma\in\mathcal E_{\rm int}\,,\sigma=K|L}\nu_K\meas(K)F_{K,\sigma}^2+
\nu_L\meas(L)F_{L,\sigma}^2+
\sum_{\sigma\in\mathcal E_{\rm ext}\,,\sigma\in\mathcal E_K}\nu_K\meas(K)F_{K,\sigma}^2\\
&=&\sum_{K\in\mathcal M}\mathbf{v}_K\cdot\sum_{\sigma\in\mathcal E_K}F_{K,\sigma}(\x_\sigma-\x_K)
+\sum_{K\in\mathcal M}\sum_{\sigma\in\mathcal E_K}\nu_K\meas(K)F_{K,\sigma}^2.
\end{eqnarray*}
Applying \refe{condconspen3}, we obtain
\begin{eqnarray}
\int_\O {\bf v}(x)\cdot \Lambda(x){\bf v}(x)\dx +\sum_{K\in\mesh}\sum_{\sigma\in\mathcal E_K}\nu_K
F_{K,\sigma}^2 \mcv &=& \int_\O f(x) u(x) \dx\label{convfort2}\\
&\le& ||f||_{L^2(\O)}||u||_{L^2(\O)}.\nonumber
\end{eqnarray}
Using Young's inequality and Lemma \ref{pdiscpen}, we deduce that,
for all $\varepsilon>0$,
\begin{eqnarray}
\alpha_0 ||\mathbf{v}||_{L^2(\O)^d}^2 + \sum_{K\in\mesh}\sum_{\sigma\in\mathcal E_K}\nu_K
\lefteqn{F_{K,\sigma}^2 \mcv \le \frac{1}{2\varepsilon}||f||_{L^2(\O)}^2
+ \varepsilon\cter{poinpen}^2||\mathbf{v}||_{L^2(\O)^d}^2}&&\nonumber\\
&&+\varepsilon \cter{poinpen}^2\sum_{K\in\mesh}\sum_{\sigma\in\mathcal E_K}\diam(K)^{2d-2}\nu_K^2F_{K,\sigma}^2
\mcv.\label{est2-je}
\end{eqnarray}
Since $\nu_K \le \nu_0 \diam(K)^{\beta}$, we have
$\nu_K \diam(K)^{2d-2} \le \nu_0 \diam(K)^{\beta + 2d-2} \le \nu_0\diam(\O)^{\beta + 2d-2}
\le \nu_0\sup(1,\diam(\O)^{\beta_0+2d-2})$ (recall that $\beta+2d-2\ge 0$).
Hence, \refe{est2-je} gives
\begin{eqnarray*}
\alpha_0 ||\mathbf{v}||_{L^2(\O)^d}^2 + \sum_{K\in\mesh}\sum_{\sigma\in\mathcal E_K}\nu_K
\lefteqn{F_{K,\sigma}^2 \mcv \le \frac{1}{2\varepsilon}||f||_{L^2(\O)}^2
+ \varepsilon\cter{poinpen}^2||\mathbf{v}||_{L^2(\O)^d}^2}&&\nonumber\\
&&+\varepsilon \nu_0 \sup(1,\diam(\O)^{\beta_0+2d-2})\cter{poinpen}^2\sum_{K\in\mesh}\sum_{\sigma\in\mathcal E_K}\nu_K F_{K,\sigma}^2
\mcv.
\end{eqnarray*}
Taking $\varepsilon=\min(\frac{\alpha_0}{2\cter{poinpen}^2},\frac{1}{2\nu_0\sup(1,\mbox{\scriptsize diam}(\O)^{\beta_0+2d-2})
\cter{poinpen}^2})$
concludes the proof of the lemma. \end{proof}

\medskip

We now prove the convergence of the approximate solution toward the weak solution
of \refe{ellgen}.

\begin{proof}{ of Theorem \ref{convpen}.}

For the simplicity of the notations, we omit the index $m$ as in the proof of
Lemma \ref{cpctpen}.
We first note that, thanks to Estimate \refe{estimpen} and since $\nu_K=\nu_0\diam(K)^\beta$,
\begin{eqnarray*}
N_2(\disc,\nu,F)^2&=&\sum_{K\in\mesh}\sum_{\sigma\in\mathcal E_K}\diam(K)^{2d-2}\nu_K^2F_{K,\sigma}^2
\mcv\\
&=&\nu_0\sum_{K\in\mesh}\sum_{\sigma\in\mathcal E_K}\diam(K)^{\beta+2d-2}\nu_KF_{K,\sigma}^2
\mcv\\
&\le& \nu_0 \size(\disc)^{\beta+2d-2} \ctel{je}
\end{eqnarray*}
where $\cter{je}$ does not depend on the discretization $\disc$ (recall that $\regul(\disc)$
is bounded). Since $\beta+2d-2>0$, this last quantity tends to $0$, and so does $N_2(\disc,\nu,F)$,
as $\size(\disc)\to 0$. Hence, still using \refe{estimpen}, we see that the assumptions of
Lemma \ref{cpctpen} are satisfied; there exists thus $\bar u\in H^1_0(\O)$
such that, up to a subsequence and as $\size(\disc)\to 0$, $\mathbf{v}\to \nabla \bar u$
weakly in $L^2(\O)^d$ and $u\to \bar u$ weakly in $L^2(\O)$
and strongly in $L^q(\O)$ for $q<2$.

We now prove that the limit function $\bar u$ is the weak
solution to \refe{ellgen}. Since any subsequence of $(u,\mathbf{v})$
has a subsequence which converges as above, and since the reasoning
we are going to make proves that any such limit of a subsequence
is the (unique) weak solution to \refe{ellgen}, this will conclude the proof, except
for the strong convergence of $\mathbf{v}$.

\medskip

Let $\varphi\in C^\infty_c(\O)$. We multiply \refe{mul3pen} by $\varphi(\x_K)$ and we sum on $K$.
Gathering by edges thanks to \refe{condconspen1}, we get
\[
\sum_{\edge\in\edges_{\rm int},\edge = K|L} F_{K,\edge} (\varphi(\x_L) -  \varphi(\x_K))= 
\sum_{K\in\mesh} \int_K \varphi(\x_K) f(x)\dx
\]
as long as $\size(\disc)$ is small enough (so that $\varphi=0$ on the control
volumes $K$ such that $\partial K\cap \partial\O\not= \emptyset$). We set, for $\edge = K|L$,
\[
\varphi(\x_L) -  \varphi(\x_K) = \frac 1 {\mcv} \int_K \grad\varphi(x)\dx
\cdot(\x_\edge - \x_K) +
\frac 1 {\meas(L)} \int_L \grad\varphi(x)\dx \cdot(\x_L - \x_\edge) + R_{KL}
\]
and we have $|R_{KL}| \le C_\varphi (\diam(K)^2+\diam(L)^2)$.
We then obtain, gathering by control volumes and using \refe{condconspen3}
(and the fact that $\varphi=0$ on the control volumes on the
boundary of $\O$),
\begin{equation}
\int_\O \Lambda_\disc {\bf v}(x)\cdot \grad\varphi(x)\dx = \int_\O f(x)
\varphi_\disc(x) \dx + \terml{convpp},
\label{estt-jd}\end{equation}
where $\Lambda_\disc$ and $\varphi_\disc$ are constant respectively equal to $\Lambda_K$ 
and $\varphi(\mathbf{x}_K)$ on each mesh $K$, and
\[
|\termr{convpp}| \le C_\varphi \sum_{\edge\in\edges_{\rm int},\edge = K|L}
|F_{K,\edge}|(\diam(K)^2+\diam(L)^2)=C_\varphi \sum_{K\in\mesh}\sum_{\sigma\in\mathcal E_K}
\diam(K)^2 |F_{K,\sigma}|.
\]
Let us estimate this term. We have
\begin{eqnarray}
|\termr{convpp}|^2 &\le& C_\varphi^2 \left(\sum_{K\in\mesh}\sum_{\edge\in\edgescv} \nu_K
F_{K,\edge}^2 \mcv\right)
\left(\sum_{K\in\mesh}\sum_{\sigma\in\mathcal E_K} \frac {\diam(K)^4} {\nu_K\mcv}\right)\nonumber\\
&\le& \ctel{truc}\sum_{K\in\mesh}\sum_{\sigma\in\mathcal E_K} \frac{\diam(K)^4}{\nu_K\mcv^2} \mcv
\label{estfinal-je}\end{eqnarray}
where, according to \refe{estimpen}, $\cter{truc}$ does not depend
on the mesh since $\regul(\disc)$ stays bounded.
But $\nu_K=\nu_0 \diam(K)^{\beta}$ and $\diam(K)^d\le \frac{{\rm regul}(\disc)}{\omega_d}\mcv$, so that
\[
\frac{\diam(K)^4}{\nu_K\mcv^2}\le\frac{\regul(\disc)^2\diam(K)^{4-\beta}}{\omega_d^2\nu_0\diam(K)^{2d}}
=\frac{\regul(\disc)^2}{\omega_d^2\nu_0}\diam(K)^{4-2d-\beta}.
\]
Since $4-2d-\beta>0$, we deduce from \refe{estfinal-je} that
\[
|\termr{convpp}|^2 \le \cter{truc}\frac{\regul(\disc)^2}{\omega_d^2\nu_0}
\size(\disc)^{4-2d-\beta} \sum_{K\in\mesh}\sum_{\sigma\in\mathcal E_K} \mcv
\le \frac{\cter{truc}\regul(\disc)^3\meas(\Omega)}{\omega_d^2\nu_0}\size(\disc)^{4-2d-\beta}
\]
and this quantity tends to $0$ as $\size(\disc)\to 0$. Hence, we can pass to the limit
in \refe{estt-jd} to see that
\[
\int_\O \Lambda \nabla\bar u(x)\cdot \grad\varphi(x)\dx = \int_\O f(x)\varphi(x) \dx\,,
\]
which proves that $\bar u$ is the weak solution to \refe{ellgen}.

\medskip

The strong convergence of $\mathbf{v}$ to $\nabla \bar u$ is a consequence of \refe{convfort2}.
{}From this equation, and defining $N(\mathbf{w})^2=\int_\O \Lambda(x)\mathbf{w}(x)\cdot\mathbf{w}(x)\dx$,
we have $N(\mathbf{v})^2\le \int_\O f(x)u(x)\dx$ and
thus
\begin{equation}
\limsup_{\mbox{\small size}(\disc)\to 0} N(\mathbf{v})^2\le 
\lim_{\mbox{\small size}(\disc)\to 0}\int_\O f(x)u(x)\dx= \int_\O f(x)\bar u(x)\dx=N(\nabla \bar u)^2
\label{rr}\end{equation}
(we use the fact that $u\to \bar u$ weakly in $L^2(\O)$ and that $\bar u$
is the weak solution to \refe{ellgen}). But $N$ is a norm on $L^2(\O)^d$,
equivalent to the usual norm and coming from the scalar
product $\langle \mathbf{w},\mathbf{z}\rangle = 
\int_\O \frac{\Lambda(x)+\Lambda(x)^T}{2}\mathbf{w}(x)\cdot \mathbf{z}(x)\,dx$;
since $\mathbf{v}\to \nabla \bar u$ weakly in $L^2(\O)^d$ as $\size(\disc)\to 0$, we therefore also
have $N(\nabla \bar u)\le \liminf_{\mbox{\small size}(\disc)\to 0} N(\mathbf{v})$.
We conclude with \refe{rr} that $N(\mathbf{v})\to N(\nabla \bar u)$ as $\size(\disc)\to 0$,
and thus that the weak convergence of $\mathbf{v}$ to $\nabla \bar u$ in $L^2(\O)^d$
is in fact strong. \end{proof}

\begin{remark} As a consequence of \refe{convfort2} and the strong convergence
of $\mathbf{v}$ to $\nabla \bar u$, we see that $\sum_{K\in\mesh}\sum_{\sigma\in\mathcal E_K}
\nu_K F_{K,\sigma}^2 \mcv\to 0$ as $\size(\disc)\to 0$. This strengthens Lemma 
\ref{exestpen} which only states that this quantity is bounded.
\end{remark}

\medskip 

To conclude this section, we prove the error estimates. Note that
these estimates could be extended, for $d\le 3$, to the case $\bar u\in H^2(\O)$
following some arguments of \cite{vignal}.

\begin{proof}{ of Theorem \ref{esterrorpen}.}

In this proof, we denote by $C_i$ (for all integer $i$)
various real numbers which can depend on
$d$, $\O$, $\bar u$, $\Lambda$ and $\theta$, but not on $\size(\disc)$.
We also denote, for all $K\in\mesh$ and $\sigma\in\mathcal E_K$,
$\bar u_K = \bar u(\x_K)$, $\bar u_\edge = \bar u(\x_\edge)$,
\[
\bar F_{K,\edge} = \int_\edge \Lambda(x) \grad \bar u(x) \cdot
\n_{K,\edge} \dfrontiere(x),
\]
\[
\bar {\bf v}_K = \frac 1 {\mcv} \Lambda_K^{-1} 
\sum_{\edge\in\edgescv} \bar F_{K,\edge} (\x_\edge -\x_K)
\]
(notice that $\Lambda_K$ is indeed invertible since, from
\refe{hyplambda}, $\Lambda_K\ge \alpha_0$).
Thanks to Lemma \ref{form-mag}, we have
\be
|\bar {\bf v}_K - \grad\bar u(x)|\le \ctel{cerr1bis} \diam(K),\quad\forall x\in K\,,\;\forall K\in\mesh,
\label{rez}\ee
which implies
\[
\bar {\bf v}_K\cdot (\x_\edge-\x_K)  = \bar u_\edge - \bar u_K + R_{K,\edge},\quad
\forall K\in\mesh,\ \forall \edge\in\edgescv,
\]
with  $|R_{K,\edge}|\le \ctel{cerr1} \diam(K)^2$ for all $K\in\mesh$
and $\edge\in\edgescv$.
Since $\bar u$ is a classical solution to \refe{ellgen}, we have
\[
- \sum_{\edge\in\edgescv} \bar F_{K,\edge} = \int_K f(x)\dx, \quad \forall K\in\mesh.
\]
Denoting, for all $K\in\mesh$ and
all $\edge\in\edgescv$,  $\widehat{u}_K = u_K - \bar u_K$, 
$\widehat{\bf v}_K = {\bf v}_K - \bar {\bf v}_K$
and $\widehat{F}_{K,\edge} = F_{K,\edge} - \bar F_{K,\edge}$, we see that
\be
- \sum_{\edge\in\edgescv} \widehat{F}_{K,\edge} = 0, \quad\forall K\in\mesh,
\label{err10}\ee
\be
\widehat{F}_{K,\sigma}+\widehat{F}_{L,\sigma}=0,\quad \forall \sigma=K|L\in \mathcal E_{\rm int},
\label{conschapeau}\ee
\be\ba
\dsp \meas(K) \Lambda_K \widehat{\bf v}_K = \sum_{\edge\in\edgescv} 
\widehat{F}_{K,\edge} (\x_\edge - \x_K),\quad \forall K\in\mesh,
\ea\label{err11}\ee
\be\ba
\widehat{\bf v}_K\cdot (\x_\edge-\x_K) + \widehat{\bf v}_L\cdot (\x_{L} - \x_\edge) 
+ \nu_K\mcv \widehat{F}_{K,\edge} +\nu_K\meas(K)\bar F_{K,\sigma} + R_{K,\edge} \\
\hspace*{5cm} - \nu_L\meas(L)\widehat{F}_{L,\edge} -\nu_L\meas(L)\bar F_{L,\edge}- R_{L,\edge}  = \widehat{u}_L - \widehat{u}_K , \\[0.1cm]
\hspace*{6cm} \forall K\in\mesh,\ \forall L\in\NN_K,\mbox{ with } \sigma=K|L,\\[0.2cm]
\widehat{\bf v}_K\cdot (\x_\edge-\x_K)  + \nu_K\mcv \widehat{F}_{K,\edge} +\nu_K\meas(K)\bar F_{K,\sigma}
+ R_{K,\edge} = - \widehat{u}_K ,\\
\hspace*{6cm}\forall K\in\mesh,\ \forall \edge\in\mathcal E_K\cap \mathcal E_{\rm ext}.
\ea\label{penerr12}\ee

We then get, multiplying \refe{err10} by $\widehat{u}_K$, \refe{penerr12} by 
$\widehat{F}_{K,\edge}$ and using \refe{conschapeau}, \refe{err11},
\begin{equation}
\sum_{K\in\mesh}\meas(K) \Lambda_K \widehat{\bf v}_K\cdot \widehat{\bf v}_K 
+ \sum_{K\in\mesh}\sum_{\edge\in\edgescv}\nu_K \mcv \widehat{F}_{K,\sigma}^2 =
\termr{penerr2}+ \termr{penerr3},
\label{dring}\end{equation}
where
\[
\terml{penerr2} = -\sum_{K\in\mesh}\sum_{\edge\in\edgescv}\nu_K\mcv\bar F_{K,\sigma} \widehat{F}_{K,\sigma},
\]
\[
\terml{penerr3} = - \sum_{K\in\mesh}\sum_{\edge\in\edgescv}R_{K,\edge}  \widehat{F}_{K,\sigma}.
\]

Using Young's inequality and the fact that $|\bar F_{K,\sigma}|\le
\ctel{nouv}\diam(K)^{d-1}$, we have
\begin{eqnarray*}
|\termr{penerr2}|&\le& \frac{1}{2}\sum_{K\in\mesh}\sum_{\edge\in\edgescv}\nu_K\mcv\bar F_{K,\sigma}^2
+\frac{1}{2}\sum_{K\in\mesh}\sum_{\edge\in\edgescv}\nu_K\mcv \widehat{F}_{K,\sigma}^2\\
&\le&\ctel{nnouv}\size(\disc)^{\beta+2d-2}
+\frac{1}{2}\sum_{K\in\mesh}\sum_{\edge\in\edgescv}\nu_K\mcv \widehat{F}_{K,\sigma}^2.
\end{eqnarray*}
Similarly, since $|R_{K,\edge}|\le \cter{cerr1} \diam(K)^2$,
\begin{eqnarray*}
|\termr{penerr3}|&\le& \frac{1}{2}\sum_{K\in\mesh}\sum_{\edge\in\edgescv}\frac{R_{K,\edge}^2}{\nu_K\meas(K)}
+\frac{1}{2}\sum_{K\in\mesh}\sum_{\edge\in\edgescv}\nu_K\meas(K)\widehat{F}_{K,\sigma}^2\\
&\le& \cter{cerr1}^2\sum_{K\in\mesh}\sum_{\edge\in\edgescv}\frac{\diam(K)^4}{\nu_K\meas(K)^2}\meas(K)
+\frac{1}{2}\sum_{K\in\mesh}\sum_{\edge\in\edgescv}\nu_K\meas(K)\widehat{F}_{K,\sigma}^2\\
&\le& \ctel{ccerr1}\size(\disc)^{4-2d-\beta}
+\frac{1}{2}\sum_{K\in\mesh}\sum_{\edge\in\edgescv}\nu_K\meas(K)\widehat{F}_{K,\sigma}^2.
\end{eqnarray*}
Gathering these two estimates in \refe{dring}, the terms involving
$\widehat{F}_{K,\sigma}$ in the left-hand side and the right-hand side compensate
and we obtain
\begin{equation}
\alpha_0||\widehat{\bf v}||_{L^2(\Omega)^d}^2
\le \ctel{ppmm}\left(\size(\disc)^{\beta+2d-2}+\size(\disc)^{4-2d-\beta}\right).
\label{tpa}\end{equation}
Estimate \refe{peneqerr0} follows, using the fact that $\size(\disc)\le 1$ and
that $||\bar\mathbf{v}-\nabla\bar u||_{L^\infty(\O)^d}\le \cter{cerr1bis}\size(\disc)$.

We now set $\widetilde{F}_{K,\edge} = \widehat{F}_{K,\edge} +\bar F_{K,\sigma}+ \frac {R_{K,\edge}} {\nu_K \mcv}=
F_{K,\edge} + \frac {R_{K,\edge}} {\nu_K \mcv}$ for all
$K\in\mesh$ and $\edge\in\edges$, and we estimate $N_2(\disc,\nu,\widetilde{F})$ the following way:
\begin{eqnarray}
N_2(\disc,\nu,\widetilde{F})^2& = &\sum_{K\in\mesh}
\sum_{\edge\in\edgescv} \diam(K)^{2d-2}\nu_K^2 \widetilde{F}_{K,\sigma}^2 \mcv\nonumber\\
&\le& \dsp 2 \sum_{K\in\mesh}
\sum_{\edge\in\edgescv} \diam(K)^{2d-2}\nu_K^2 F_{K,\sigma}^2 \mcv \nonumber\\ 
&&+2 \sum_{K\in\mesh} \sum_{\edge\in\edgescv} \diam(K)^{2d-2} \nu_K^2 \mcv \frac 
{\cter{cerr1}^2 \diam(K)^4} {(\nu_K \mcv)^2}\nonumber\\
&\le& \ctel{penerr6} (\size(\disc)^{\beta + 2 d - 2} +\size(\disc)^2)
\label{penestilF}
\end{eqnarray}
(we have used \refe{estimpen}).
Since  \refe{penerr12} implies that $(\widehat{u},\widehat{\bf v},\widetilde{F})\in L_\nu(\disc)$,
Lemma \ref{pdiscpen} gives
\[
\Vert \widehat{u}\Vert_{L^2(\O)} \le 
\cter{poinpen}\left(\Vert \widehat{\bf v}\Vert_{L^2(\O)^d}+
N_2(\disc,\nu,\widetilde{F})\right)
\]
and \refe{peneqerr1} follows from \refe{tpa}, \refe{penestilF}
and an easy estimate between $\bar u_K$ and the values of $\bar u$
on $K$. \end{proof}

\section{Implementation}\label{numimp}

We present the practical implementation in the case
where $\Lambda(x)$ is symmetric for a.e. $x\in\O$, though it is
valid for any $\Lambda$ (notice that, in the physical problems given
in the introduction of this paper,
the diffusion tensor is always symmetric).

\subsection{Resolution procedure}

The size of System (\refe{condconspen2},\refe{condconspen1},\refe{condconspen3},\refe{mul3pen}) is equal
to $(d+1)\mbox{Card}(\mesh) + 2\mbox{Card}(\edgesint) + \mbox{Card}(\edgesext)$. However, it is possible
to proceed to an algebraic elimination which leads to a 
symmetric positive definite sparse linear
system with $\mbox{Card}(\edgesint)$ unknowns, following
the same principles as in the hybrid resolution of a mixed finite
element problem (see for example \cite{RobT}). Indeed, 
for all $(u,{\bf v},F)$ such that \refe{condconspen2}
and \refe{condconspen3} hold, we define
$(u_\edge)_{\edge\in\edges}$ by
\[
{\bf v}_K\cdot (\x_\edge-\x_K) + \nu_K F_{K,\sigma} \mcv =
u_\sigma- u_K,\quad
\forall K\in\mesh,\ \forall \edge\in\edgescv.
\]
We thus have that $u_\edge = 0$ for all $\edge\in\mathcal E_K\cap \mathcal E_{\rm ext}$.
We can then express $({\bf v},F)$ as a function of 
$(u_\edge)_{\edge\in\edges}$ and of $u$, since we have
\[\ba \dsp
\frac 1 {\mcv}\sum_{\edge'\in\edgescv} F_{K,\edge'}  \Lambda_K^{-1}(\x_{\edge'} -
\x_K)\cdot (\x_\edge-\x_K) + \nu_K F_{K,\sigma} \mcv = u_\sigma- u_K,\\
\hspace*{7cm}\forall K\in\mesh,\ \forall \edge\in\edgescv,
\ea\]
which is, for all $K\in\mesh$, an invertible linear system with unknown
$(F_{K,\edge})_{\edge\in\edgescv}$, under the form 
$B_K(F_{K,\edge})_{\edge\in\edgescv} = (u_\edge- u_K)_{\edge\in\edgescv}$
where $B_K$ is a symmetric positive definite matrix (thanks
to the condition $\nu_K >0$).
We can then write 
\be
F_{K,\edge} = \sum_{\edge'\in\edgescv} (B_K^{-1})_{\edge\edge'}(u_{\edge'}- u_K),
\quad \forall K\in\mesh,\ \forall \edge\in\edgescv.
\label{calfks}\ee
We then obtain from \refe{mul3pen}, 
denoting $b_{K,\edge'} = \sum_{\edge\in\edgescv} (B_K^{-1})_{\edge\edge'}$ 
and $b_K = \sum_{\edge'\in\edgescv}b_{K,\edge'}$, that $u_K$ satisfies the relation
\be
- \sum_{\edge'\in\edgescv}b_{K,\edge'} u_{\edge'} + b_{K}u_K = \int_K f(x)\dx.
\label{caluk}\ee
We have $(b_{K,\edge'})_{\edge'\in\edgescv} = B_K^{-1} (1)_{\edge'\in\edgescv}$
and therefore we get  $b_K = (1)_{\edge'\in\edgescv} \cdot  B_K^{-1} 
(1)_{\edge'\in\edgescv}> 0$ since $B_K^{-1}$ is symmetric positive definite.
Reporting the previous linear relations in \refe{condconspen1}, we find
\be\ba\dsp
\sum_{\edge'\in\edgescv}
\left( (B_K^{-1})_{\edge\edge'}- \frac {b_{K,\edge}b_{K,\edge'}} {b_K} \right) u_{\edge'} +
\sum_{\edge'\in\edges_L}\left( (B_L^{-1})_{\edge\edge'} - \frac {b_{L,\edge}b_{L,\edge'}} {b_L} \right) 
u_{\edge'} = \\[0.5cm]
\qquad \dsp \frac {b_{K,\edge}} {b_K} \int_{K} f(x)\dx + \frac {b_{L,\edge}} {b_L} \int_{L} f(x)\dx,
\quad \forall \edge=K|L \in\edgesint,
\ea\label{aresoudre}\ee
which is a symmetric linear system, whose unknowns are 
$(u_\edge)_{\edge\in\edgesint}$. Let us show that its matrix $M$
is positive. We can write, for all 
family of real numbers $(u_\edge)_{\edge\in\edgesint}$,
\[
(u_\edge)_{\edge\in\edgesint}\cdot M\ (u_\edge)_{\edge\in\edgesint}
= \sum_{K\in\mesh} 
\left(
\sum_{\edge\in\edgescv}\sum_{\edge'\in\edgescv}
 (B_K^{-1})_{\edge\edge'}u_\edge u_{\edge'}  - \frac {(\sum_{\edge\in\edgescv}b_{K,\edge} u_\edge)^2} {b_K}\right).
\]
Thanks to the fact that $B_K^{-1}$ is symmetric positive definite, we get, using the 
Cauchy-Schwarz inequality,
\[
\left((1)_{\edge\in\edgescv} \cdot B_K^{-1} (u_\edge)_{\edge\in\edgescv}
\right)^2 \le \left((1)_{\edge\in\edgescv} \cdot B_K^{-1} (1)_{\edge\in\edgescv}
\right)\left((u_\edge)_{\edge\in\edgescv} \cdot B_K^{-1} (u_\edge)_{\edge\in\edgescv}
\right),
\]
which is exactly 
\[
\left(\sum_{\edge\in\edgescv}b_{K,\edge} u_\edge\right)^2 \le b_K 
\sum_{\edge\in\edgescv}\sum_{\edge'\in\edgescv}
 (B_K^{-1})_{\edge\edge'}u_\edge u_{\edge'}.
\]
In order to show that $M$ is definite, we simply remark that the preceding
reasoning shows that the systems (\refe{condconspen2},\refe{condconspen1},\refe{condconspen3},\refe{mul3pen})
and \refe{aresoudre} are equivalent. Hence, since
(\refe{condconspen2},\refe{condconspen1},\refe{condconspen3},\refe{mul3pen}) has
a unique solution, so must \refe{aresoudre}, which means that $M$ is invertible.

\medskip

Hence, we can first solve $(u_\edge)_{\edge\in\edgesint}$ from \refe{aresoudre},
and then compute $(u,F)$ thanks to relations \refe{caluk} and \refe{calfks}
and finally $\mathbf{v}$ by \refe{condconspen3}.

\subsection{Numerical results}\label{secnum}

Taking $\nu_K = 0$ for all $K\in\mesh$, we could prove in the symmetric
case, via a minimization technique, that
there exists at least one $(u,{\bf v},F)\in L_\nu(\disc)$ solution of 
(\refe{condconspen2},\refe{condconspen1},\refe{condconspen3},\refe{mul3pen}).
In this case, $(u,{\bf v})$ is unique, but
this is no longer true for $F$ in the general case (see however 
section \ref{sixdeux}). Within such a choice,
the proof of convergence of $(u,{\bf v})$ to the continuous
solution remains an open problem.
Nevertheless,  this gives an indication that very small values
of $(\nu_K)_{K\in\mesh}$ can be considered. Hence
we take $\nu_K = 10^{-9}/ \mcv$ in all the following computations.
The inversion of matrices $B_K$ arising in \refe{calfks}
and the solving of System \refe{aresoudre} are then realized using
direct methods.

\subsubsection{Case of a homogeneous isotropic problem}

We consider here the case $d=2$, $\O = (0,1)\times(0,1)$,
$\Lambda = {\rm I}_d$ and $\bar u(x) = x_1(1- x_1) x_2(1- x_2)$
for all $x =  (x_1, x_2)\in\O$.

\begin{figure}
\begin{center}\begin{tabular}{cc}
\includegraphics[width=0.25 \linewidth]{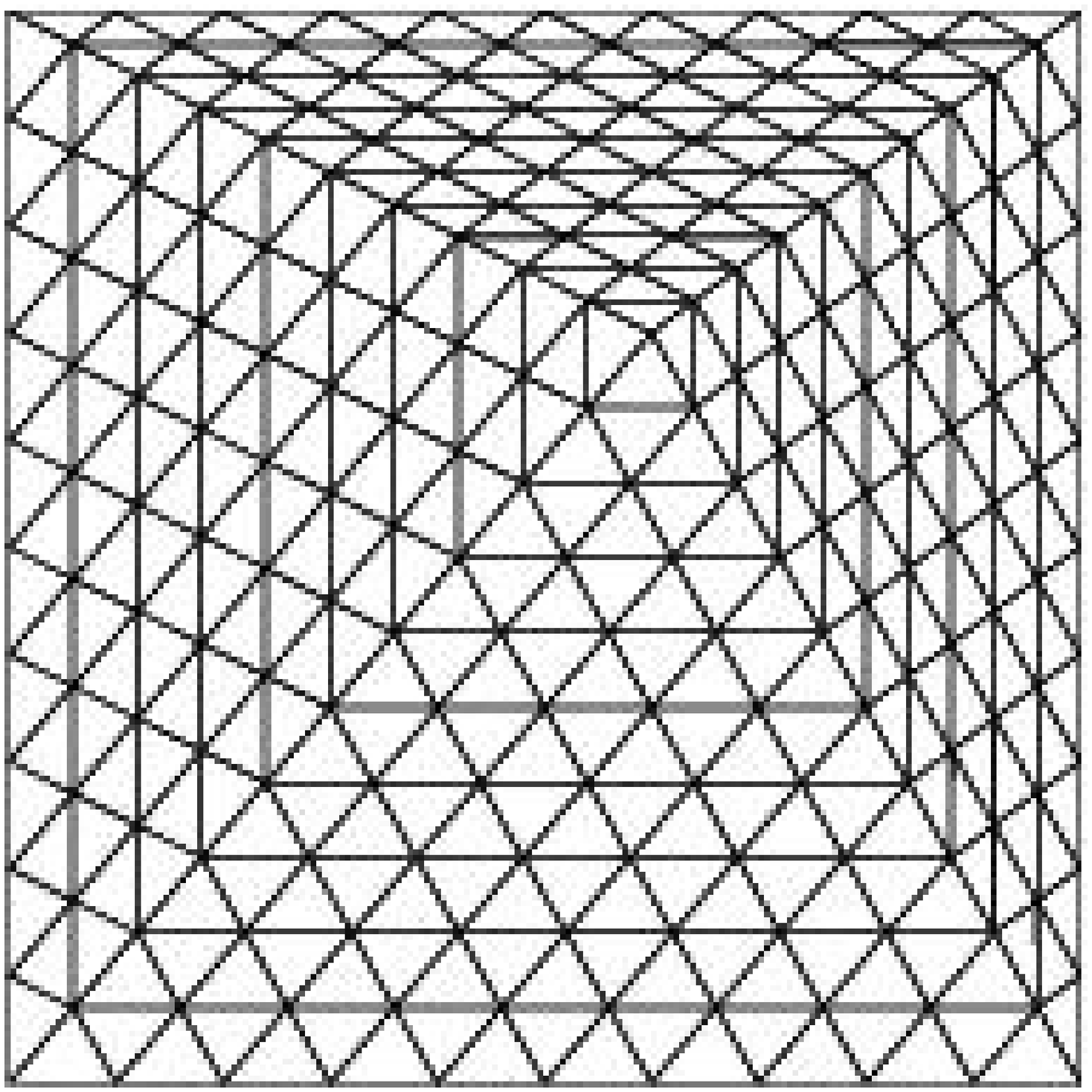}&
\includegraphics[width=0.25 \linewidth]{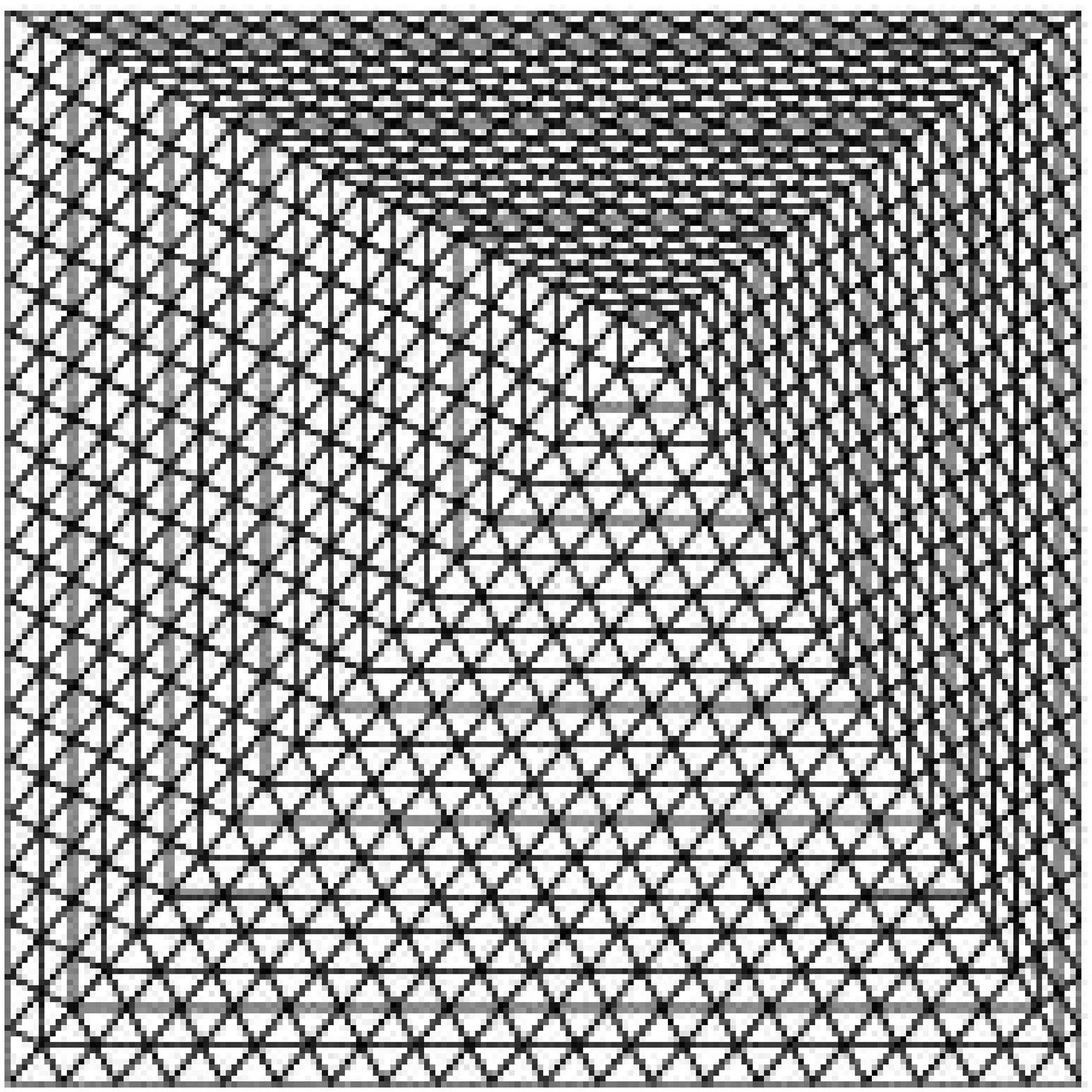}\\
grid $\disc_{t1}$ & grid $\disc_{t2}$ \\ \\
\includegraphics[width=0.25 \linewidth]{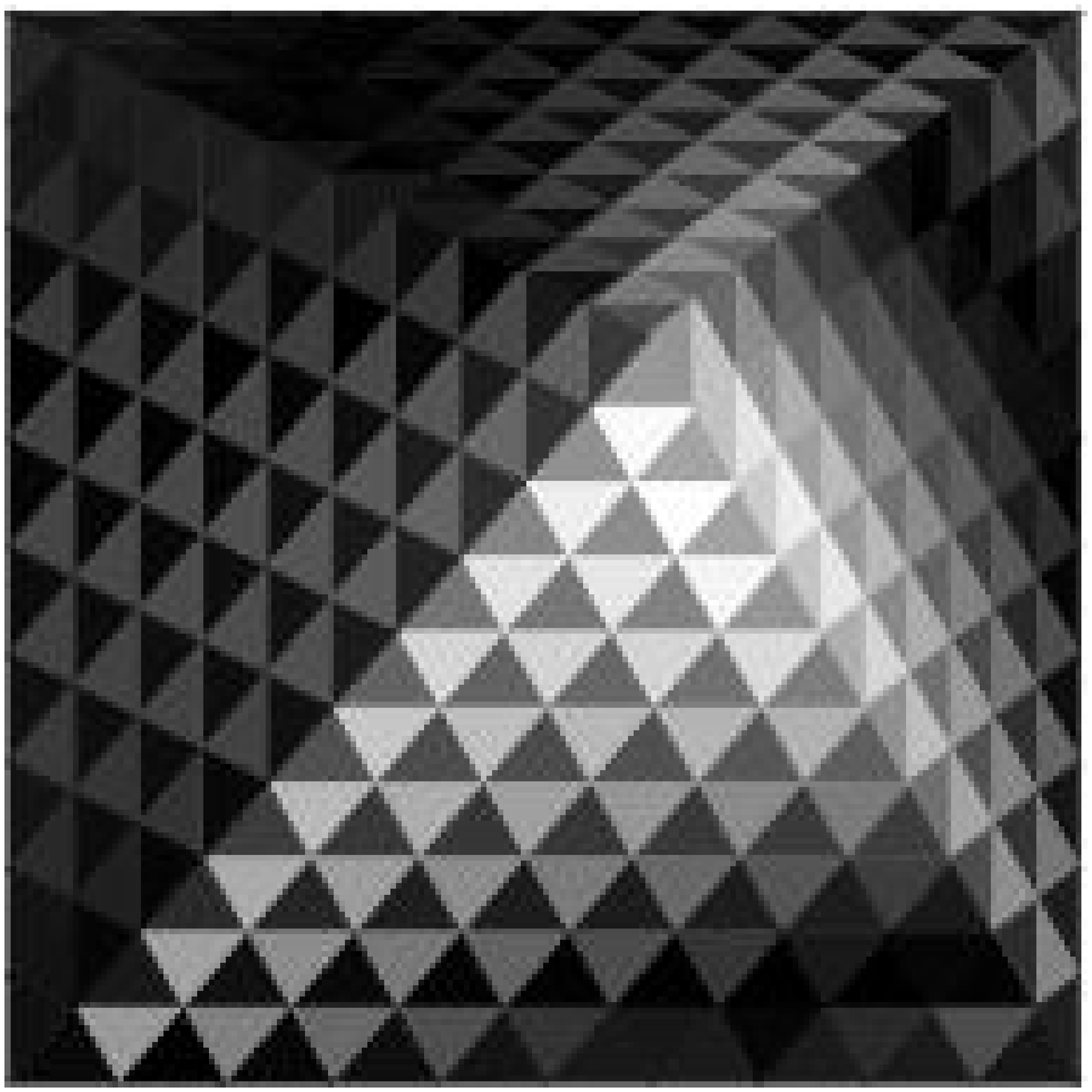}&
\includegraphics[width=0.25 \linewidth]{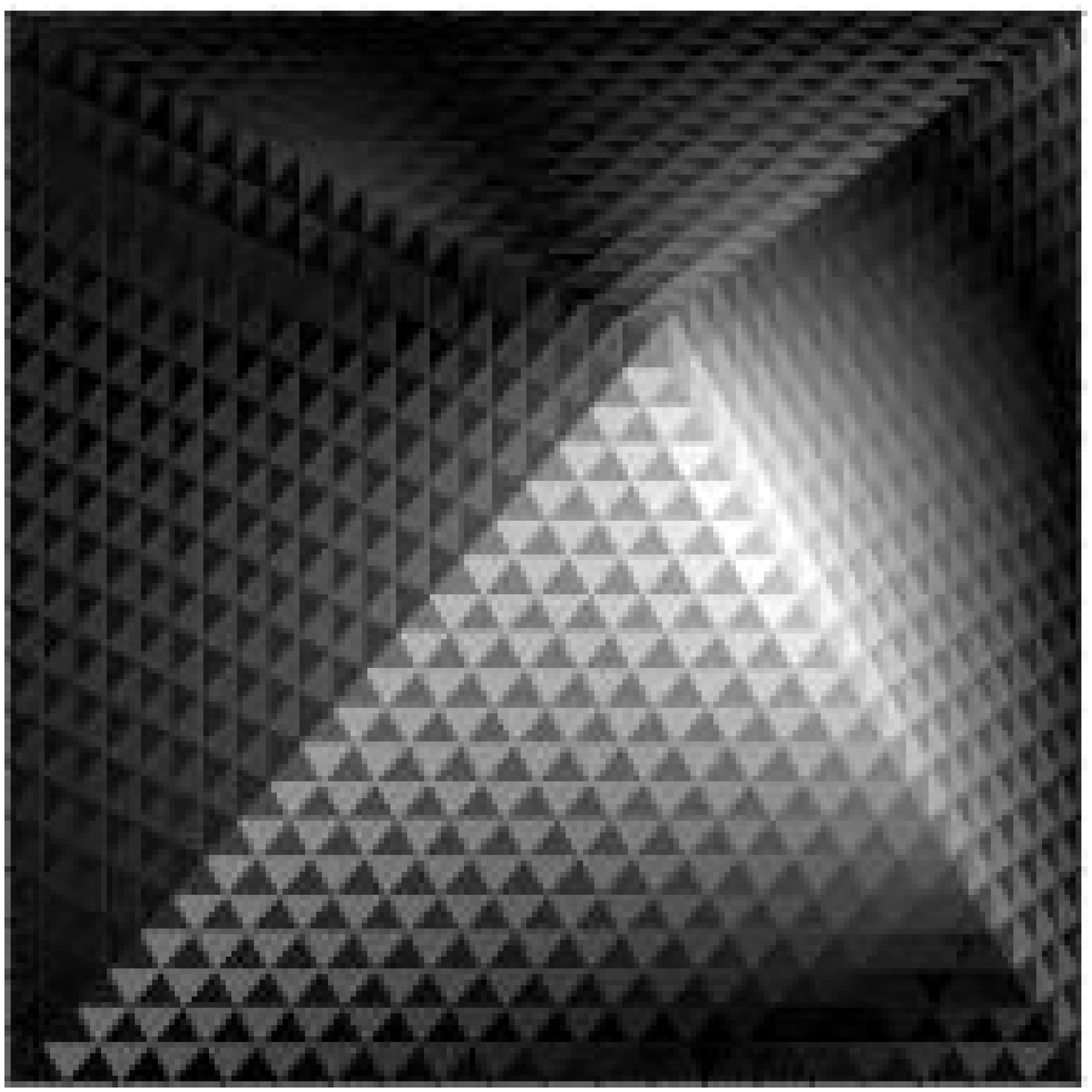}\\
error on $\disc_{t1}$ & error on $\disc_{t2}$ \\
black=0, white = $2.2\ 10^{-2}$& black=0, white = $8.9\ 10^{-3}$
\end{tabular}
\end{center}\caption{Discretizations and error  $e_\disc$ on grids
$\disc_{t1}$ and $\disc_{t2}$.}
\label{gridt}
\end{figure}

We first present in Figure \ref{gridt} two different triangular
discretizations  $\disc_{t1}$ and  $\disc_{t2}$ used for
the computation of an approximate solution for the problem.
We also show in Figure \ref{gridt} the error  $e_\disc$, defined by
\[
e_K = \frac {|u_K - \bar u(\x_K)|} {\Vert \bar u\Vert_{L^\infty(\O)}},
\quad \forall K\in\mesh,
\]
using  discretizations $\disc_{t1}$ and  $\disc_{t2}$. Note that
these discretizations do not respect the Delaunay condition on a sub-domain of $\O$, and that the
4-point finite volume scheme (see \cite{book}) cannot be used on these grids.
The grids $\disc_{t2}$ and  $\disc_{t3}$ (which is not represented here) 
have been obtained from  $\disc_{t1}$ (containing 400 control volumes) by 
the respective divisions by 2 and 4 of each edge (there are 1600 control volumes in
$\disc_{t2}$  and 6400 in  $\disc_{t3}$). For all these discretizations, 
the points $\x_K$ have been located at the center of gravity of
the control volumes.
The errors in $L^2$ norms obtained with these grids 
are given in the  following table.
\begin{center}
\begin{tabular}{| c | c | c |}
\hline
 & $\Vert u - \bar u\Vert_{L^2(\O)}$  &  $\Vert {\bf v} - \grad \bar u\Vert_{L^2(\O)^d}$ \\
\hline
$\disc_{t1}$       &  $5.1\ 10^{-4}$    &  $1.8\ 10^{-2}$    \\
\hline
$\disc_{t2}$       &  $1.9\ 10^{-4}$    &  $9.0\ 10^{-3}$    \\
\hline
$\disc_{t3}$       &  $8.2\ 10^{-5}$    &  $4.5\ 10^{-3}$    \\
\hline
order of convergence     &   $\ge 1$  &   $1$   \\
\hline
\end{tabular}
\end{center}

We observe that the  numerical orders of convergence for
$\Vert u - \bar u\Vert_{L^2(\O)}$ and 
$\Vert {\bf v} - \grad \bar u\Vert_{L^2(\O)^d}$  seem to be equal to 1, 
and therefore
no super-convergence property can reasonably be expected in this case.

\begin{figure}
\begin{center}\begin{tabular}{cc}
\includegraphics[width=0.25 \linewidth]{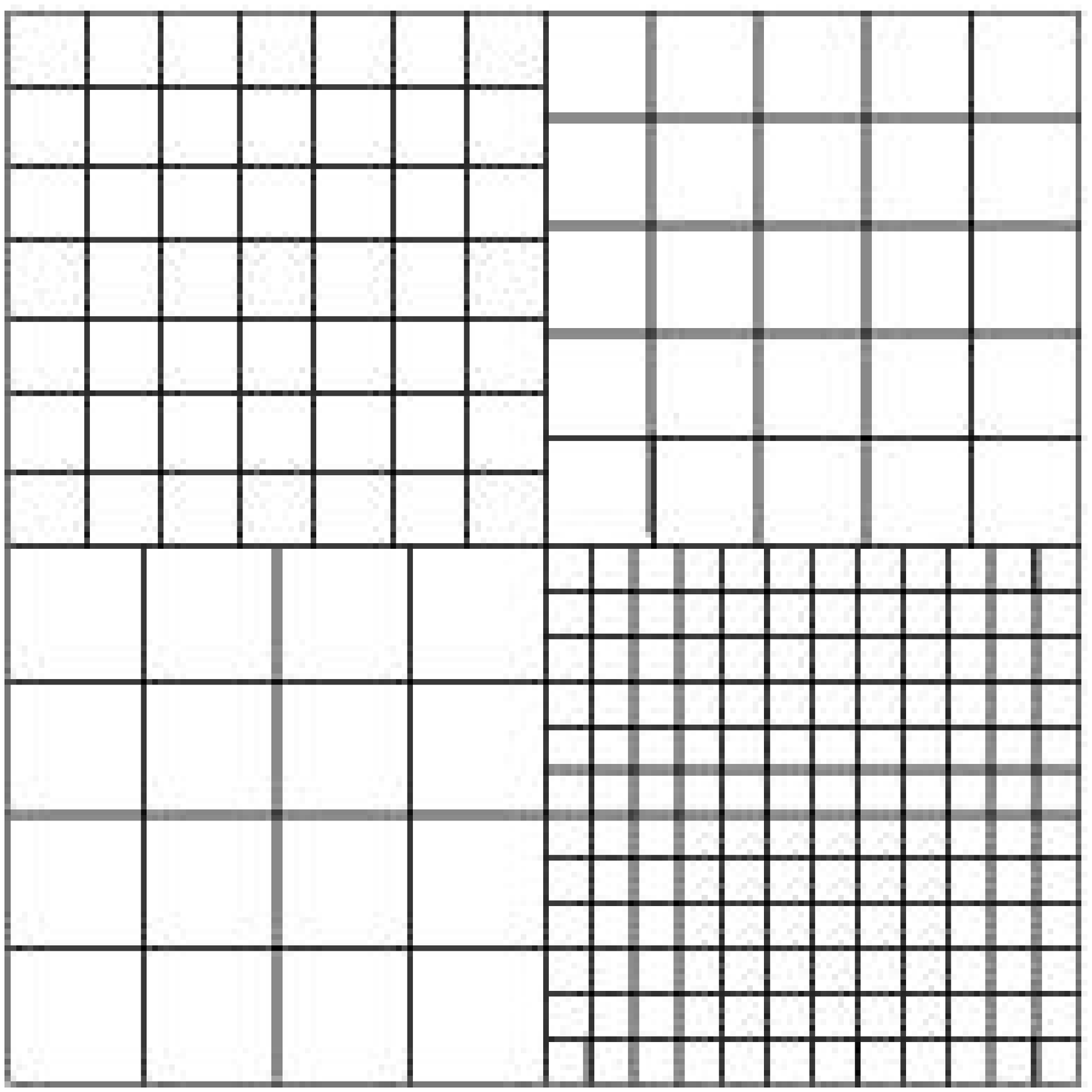}&
\includegraphics[width=0.25 \linewidth]{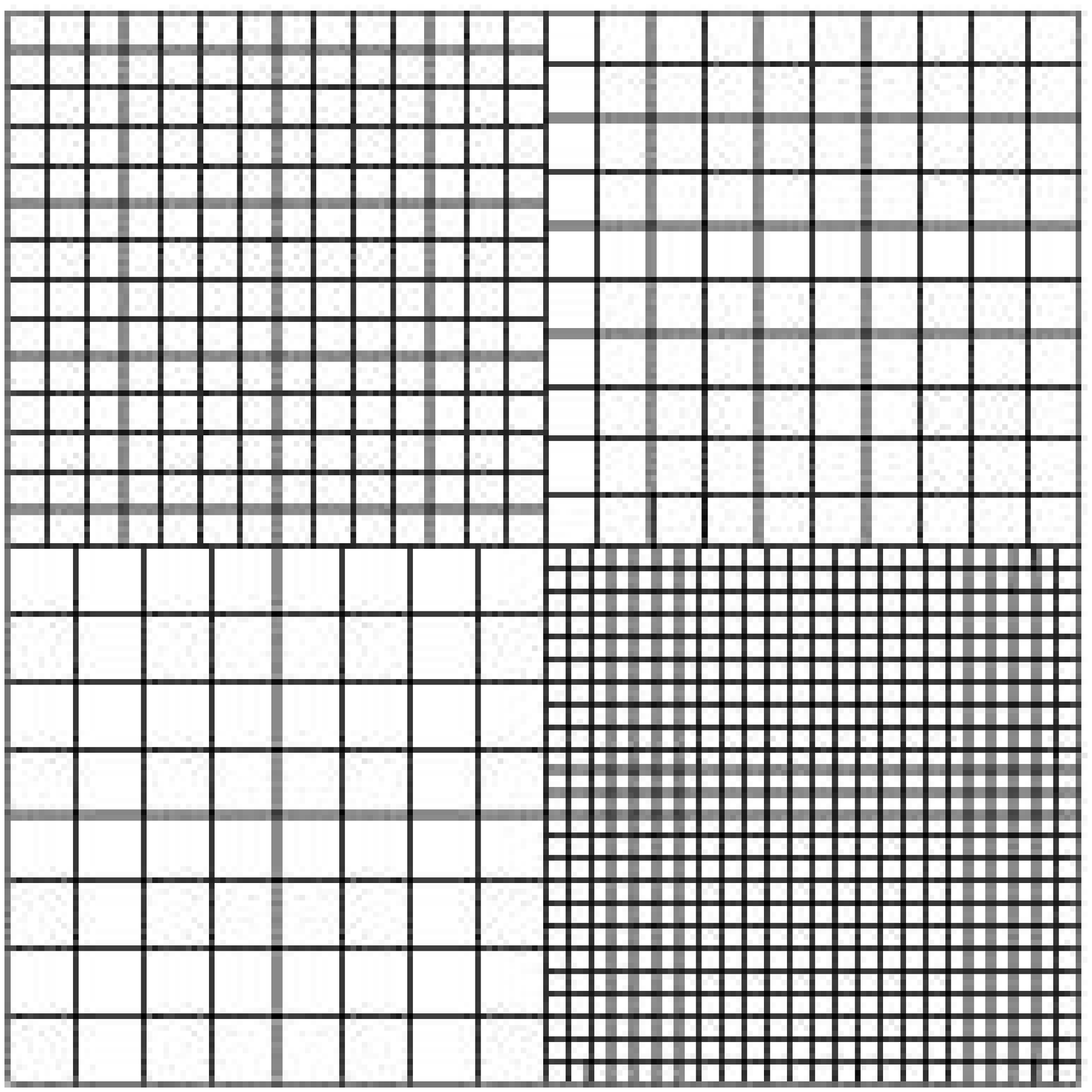}\\
grid $\disc_{q1}$ & grid $\disc_{q2}$ \\ \\
\includegraphics[width=0.25 \linewidth]{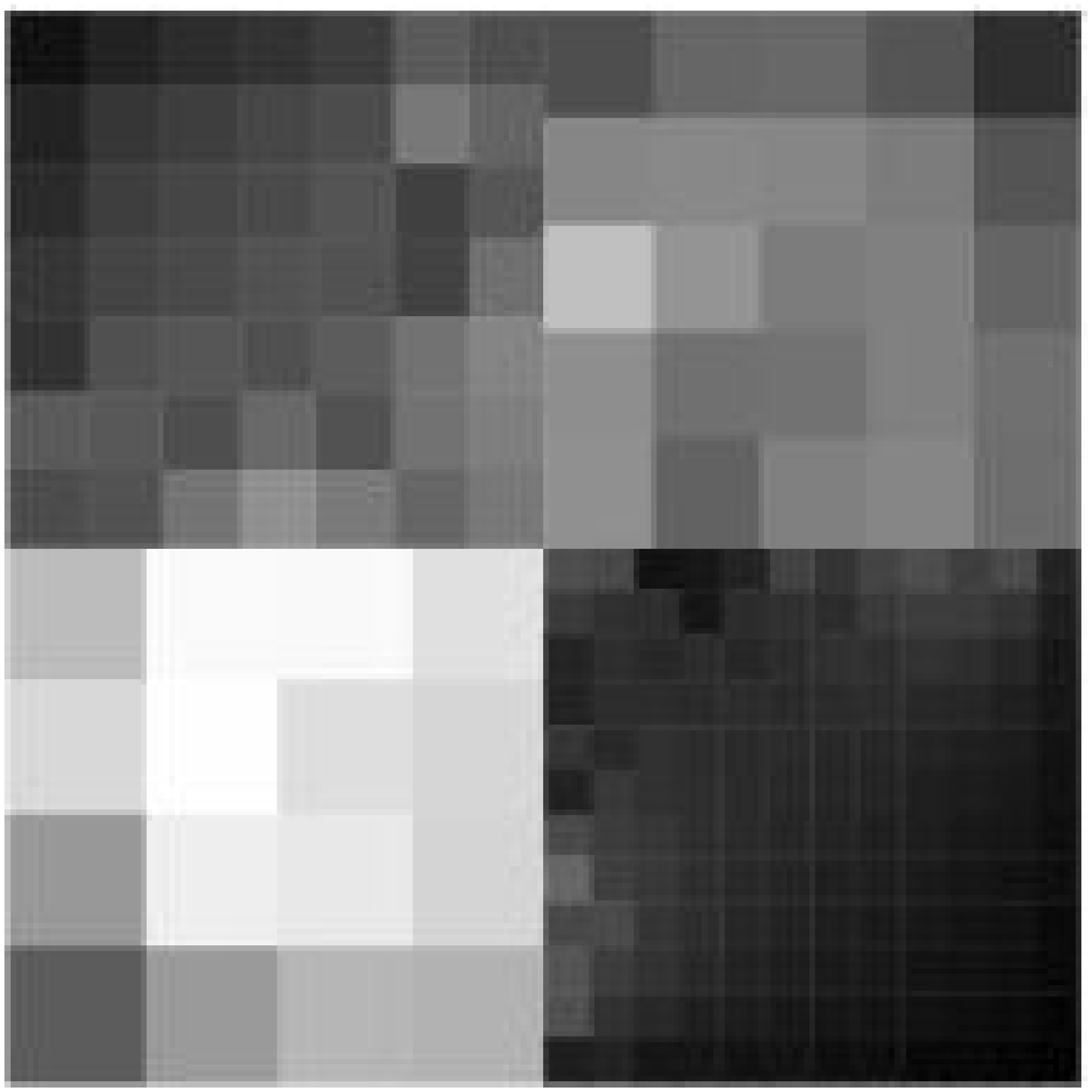}&
\includegraphics[width=0.25 \linewidth]{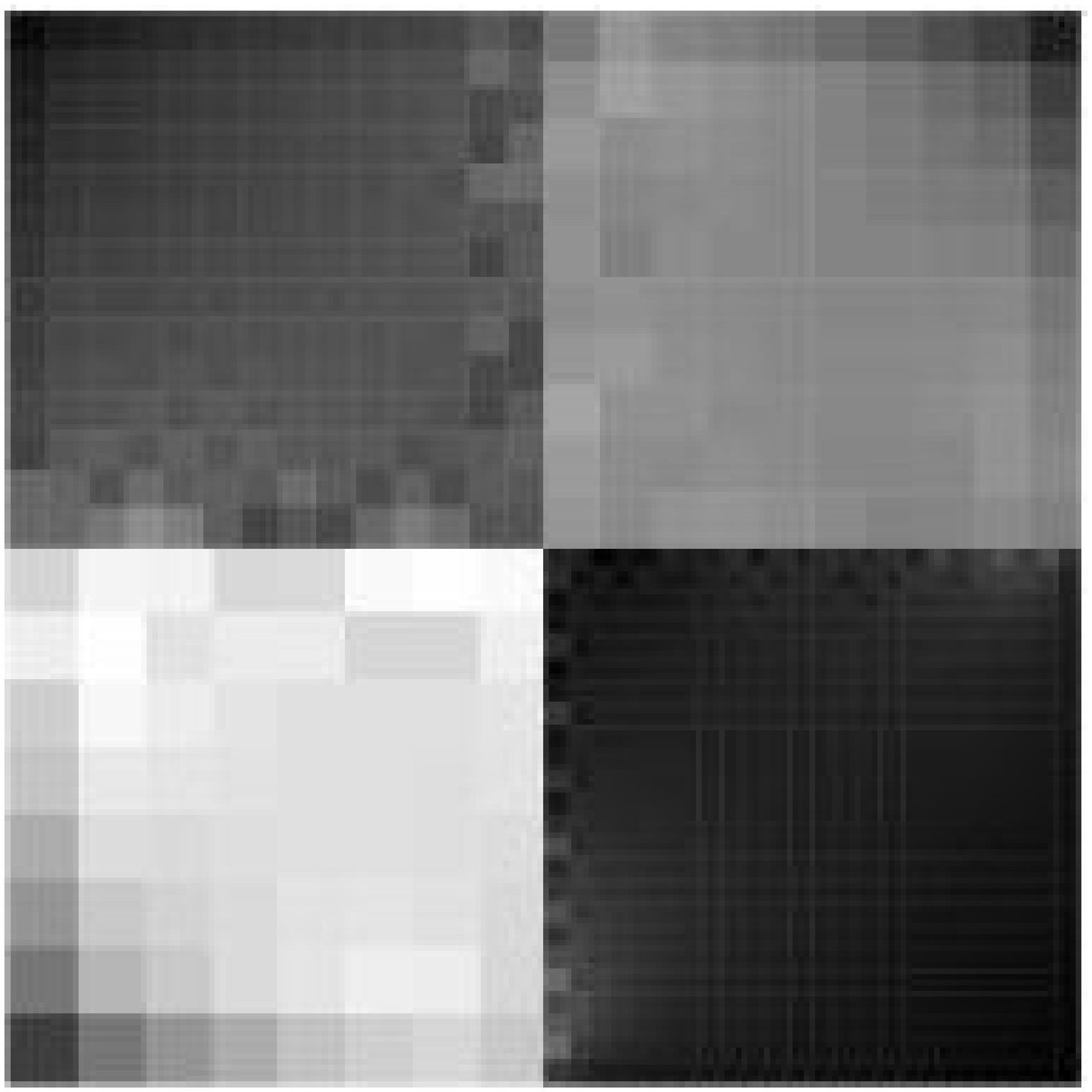}\\
error on $\disc_{q1}$ & error on $\disc_{q2}$ \\
black=0, white = $2.7\ 10^{-2}$& black=0, white = $5.3\ 10^{-3}$
\end{tabular}
\end{center}\caption{Discretizations and error  $e_\disc$ on grids
$\disc_{q1}$ and $\disc_{q2}$.}
\label{grids}
\end{figure}

We then present in Figure \ref{grids} discretizations 
$\disc_{q1}$ and $\disc_{q2}$ and error $e_\disc$ using these grids.
Such grids could be obtained using a refinement procedure: for example,
in the case of coupled systems, the grid might have been refined
in order to improve the convergence on another equation (thanks to some
\emph{a posteriori} estimates maybe) and
must then be used to solve \refe{ellgen} which is the second part of the system.
The grid  $\disc_{q2}$ has been obtained from $\disc_{q1}$
by a uniform division of each edge by 2, and similarly  
$\disc_{q3}$ (not represented here)
has been obtained from $\disc_{q2}$ in the same way.
The respective errors in $L^2$ norms obtained with these grids 
are given in the  following table.

\begin{center}
\begin{tabular}{| c | c | c |}
\hline
 & $\Vert u - \bar u\Vert_{L^2(\O)}$  &  $\Vert {\bf v} - \grad \bar u\Vert_{L^2(\O)^d}$ \\
\hline
$\disc_{q1}$       &  $8.7\ 10^{-4}$    &  $5.8\ 10^{-3}$    \\
\hline
$\disc_{q2}$       &  $1.7\ 10^{-4}$    &  $1.3\ 10^{-3}$    \\
\hline
$\disc_{q3}$       &  $3.9\ 10^{-5}$    &  $4.0\ 10^{-4}$    \\
\hline
order of convergence     &   $\ge 2$   &   $\ge 1$  \\
\hline
\end{tabular}
\end{center}

We then observe that the  numerical order convergence is better than 2 for
$\Vert u - \bar u\Vert_{L^2(\O)}$, which corresponds to a case of
a mainly structured grid (there is no significant 
additional error located at the internal
boundaries between the differently gridded subdomains,
see Figure \ref{grids}).

\begin{figure}
\begin{center}\begin{tabular}{cc}
\includegraphics[width=0.25 \linewidth]{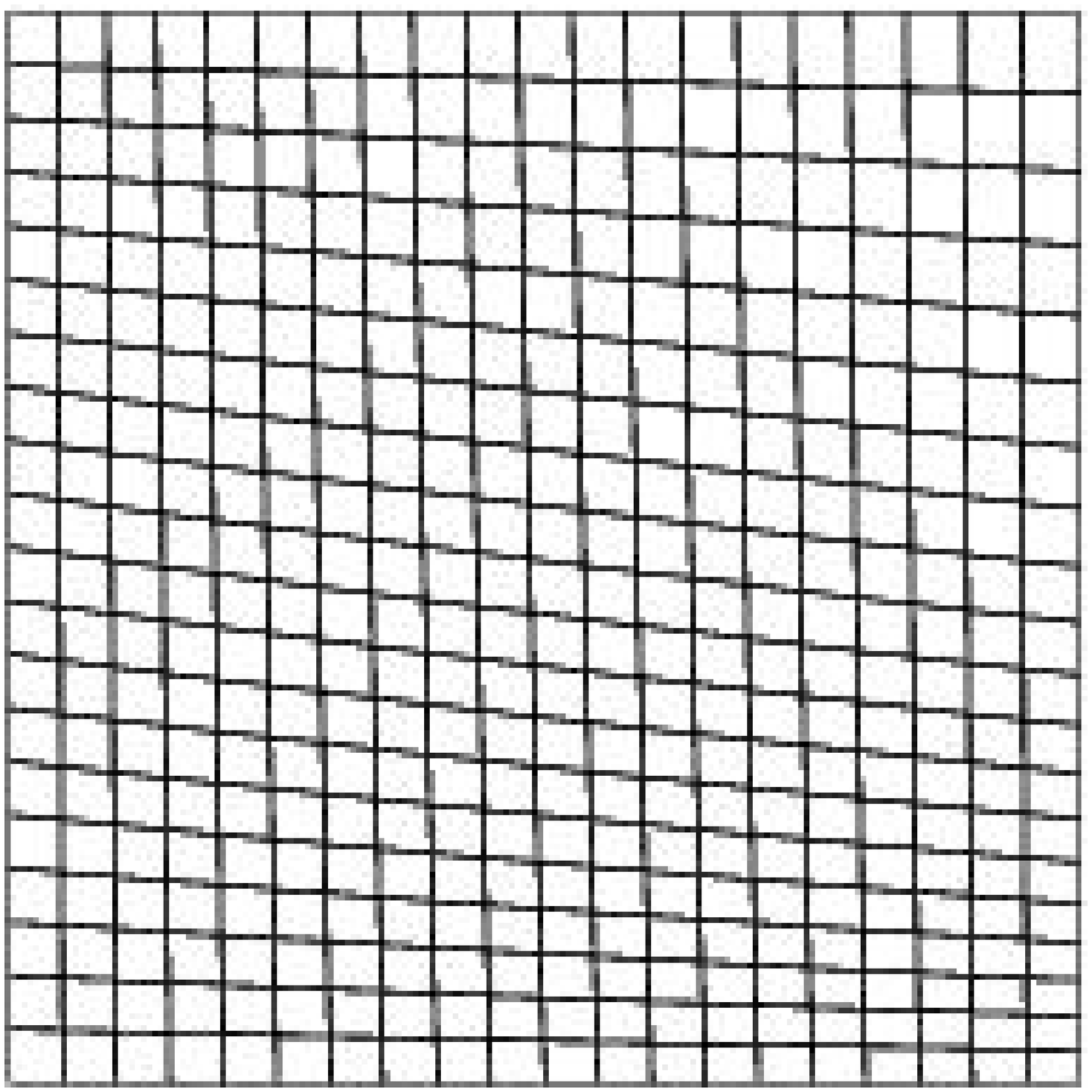}&
\includegraphics[width=0.25 \linewidth]{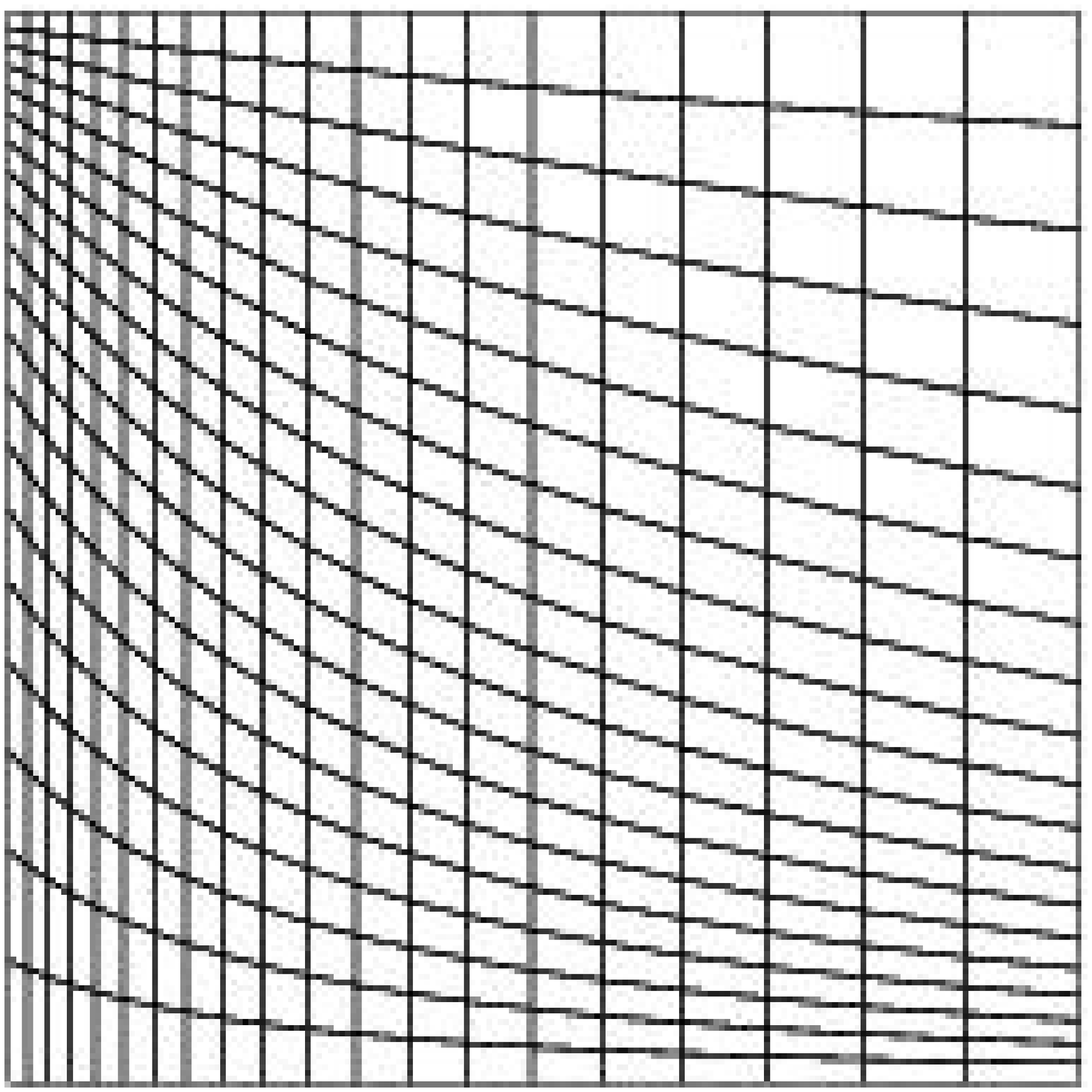}\\
grid $\disc_\flat$ & grid $\disc_\sharp$ \\ \\
\includegraphics[width=0.25 \linewidth]{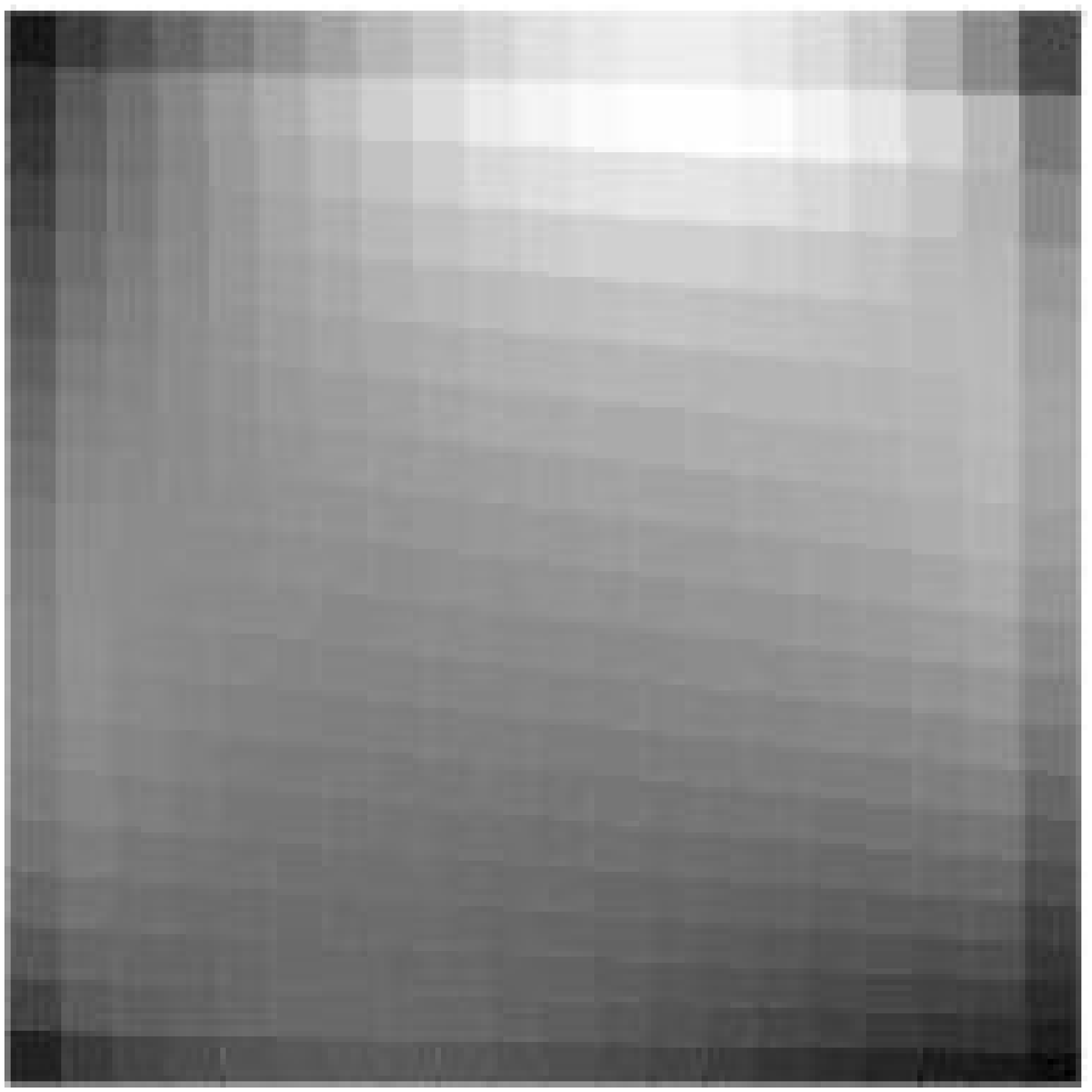}&
\includegraphics[width=0.25 \linewidth]{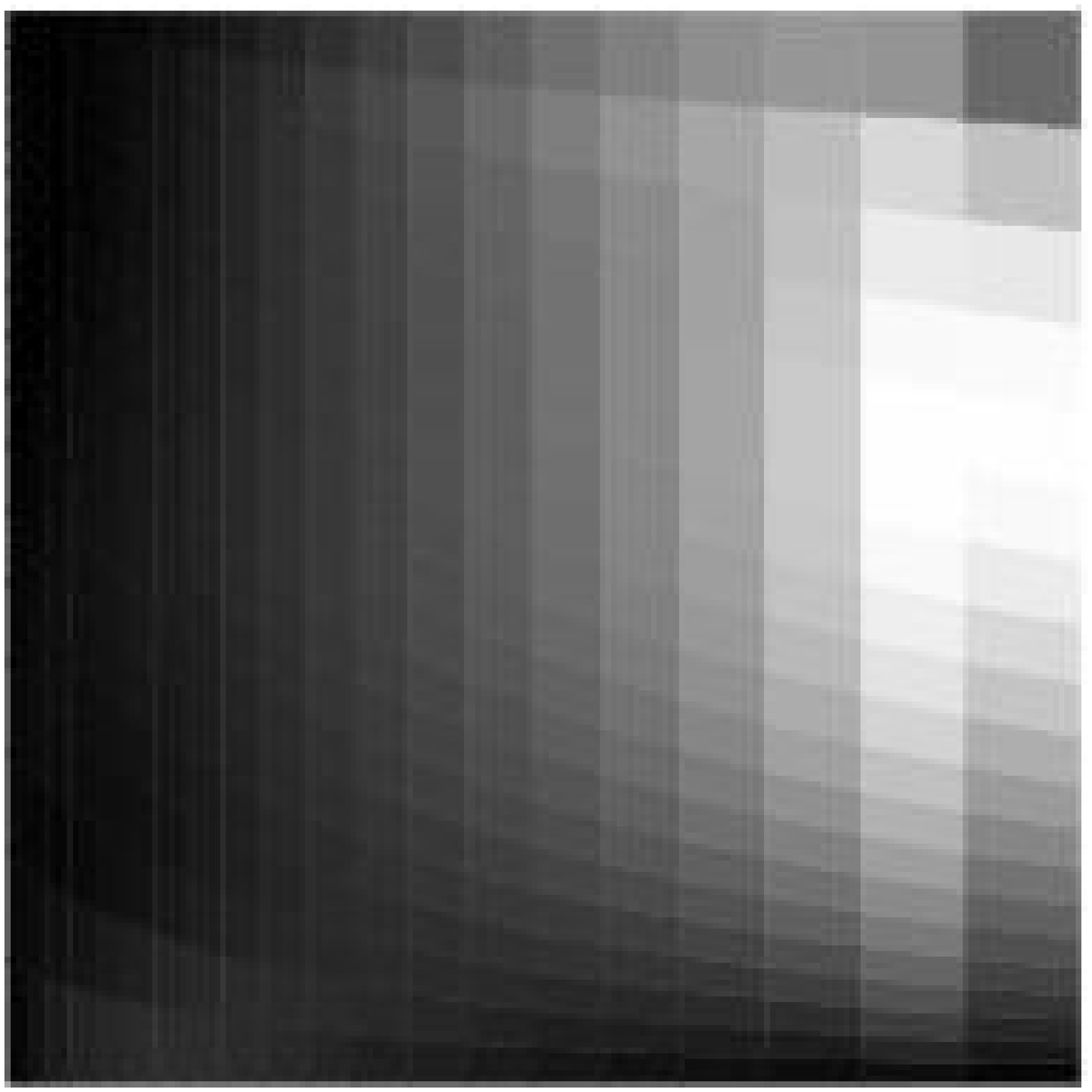} \\
error on $\disc_\flat$ & error on $\disc_\sharp$ \\
black=0, white = $5.4\ 10^{-3}$ & black=0, white = $1.5\ 10^{-2}$
\end{tabular}
\end{center}
\caption{Discretizations and error  $e_\disc$ on grids
$\disc_\flat$ and $\disc_\sharp$.}
\label{grid9}
\end{figure}

Finally, in Figure \ref{grid9}, we represent grids 
$\disc_\flat$ and $\disc_\sharp$ and the error $e_\disc$ thus 
obtained. These meshes (which have the same number of control volumes) could correspond to the case of moving meshes 
(for example, due to a phenomenon of compaction, see \cite{faille}). 
The respective errors in $L^2$ norms obtained
with these grids are given in the  following table.

\begin{center}
\begin{tabular}{| c | c | c |}
\hline
 & $\Vert u - \bar u\Vert_{L^2(\O)}$  &  $\Vert {\bf v} - \grad \bar u\Vert_{L^2(\O)^d}$ \\
\hline
$\disc_\flat$       &  $2.0\ 10^{-4}$    &  $6.7\ 10^{-4}$    \\
\hline
$\disc_\sharp$       &  $4.6\ 10^{-4}$    &  $1.8\ 10^{-3}$    \\
\hline
\end{tabular}
\end{center}

We observe that the error is mainly connected to the size of 
the control volumes, and maybe to some effect of 
loss of regularity of the mesh.

\subsubsection{Case of a heterogeneous anisotropic problem}

Let us now give some numerical results in a highly heterogeneous
and anisotropic case, inspired by \cite{lepot}. 
With $\O=(0,1)\times (0,1)$, let us
define $\bar x = (-0.1,-0.1)$ and $\eps = 10^{-4}$, and let us set
 $$\Lambda(x) = \left(\begin{array}{cc}
\dsp (x_2 -\bar x_2)^2 + \eps (x_1 -\bar x_1)^2 & 
-(1-\eps) (x_1 -\bar x_1)  (x_2 -\bar x_2)\\  & \\
\dsp
-(1-\eps) (x_1 -\bar x_1)  (x_2 -\bar x_2) & 
 (x_1 -\bar x_1)^2 + \eps (x_2 -\bar x_2)^2  \end{array} \right),\ \forall x\in\O.$$
The eigenvalues of $\Lambda(x)$ are 
equal to ${\underline \lambda}(x) =\eps \vert x - \bar x\vert^2$
and  ${\overline \lambda}(x) =\vert x - \bar x\vert^2$: the 
anisotropy ratio is therefore $1/\eps=10^4$ in  the whole domain.
Note that, thanks to the choice
$\bar x = (-0.1,-0.1)$, we have 
$\inf_{x\in\O} {\underline \lambda}(x) = \vert\bar x\vert^2\eps = 0.02 \eps$ and
$\sup_{x\in\O} {\underline \lambda}(x) = \vert \widehat x - \bar x\vert^2\eps= 2.42\eps$ 
with $\widehat x = (1,1)$. Therefore
${\underline \lambda}(x)/{\underline \lambda}(y)$ and
${\overline \lambda}(x)/{\overline \lambda}(y)$ are in the range $(1/121,121)$ for all $x,y\in \O$
(note that in \cite{lepot}, these ratios are in the range $(0,+\infty)$ since
the author takes $\bar x = (0,0)$, but then \refe{hyplambda} does not hold).
Since the directions of anisotropy are not 
constant, one cannot solve this problem
by a classical finite volume method on a tilted rectangular
mesh. 
 We assume that the solution of Problem \refe{ellgen}
is given by $\bu(x) = \sin(\pi x_1)\sin(\pi x_2)$;
in this case, $\Vert \bu\Vert_{L^2(\O)} = 1/2$, and 
the function $f$ satisfies:
 $$ \begin{array}{rcl}
f(x)  & =  & \pi^2 (1+\eps) \sin(\pi x_1)\sin(\pi x_2) \vert x - \bar x\vert^2 \\
 & &+   \pi (1 - 3\eps)\cos(\pi x_1)\sin(\pi x_2) (x_1 - \bar x_1) \\
 & & + \pi (1 - 3\eps)\sin(\pi x_1)\cos(\pi x_2) (x_2 - \bar x_2)\\
 & &+ 2\pi^2 (1-\eps) \cos(\pi x_1)\cos(\pi x_2) (x_1 - \bar x_1)(x_2 - \bar x_2),\ 
\forall x\in\O.
\end{array}$$

We then compare on this problem the numerical solution
given by the mixed finite volume scheme (denoted
by MFV below), with that obtained using the low degree mixed finite element scheme (denoted
by MFE below) in the case of triangles or rectangles. 
We compute the solutions with both schemes on the following grids: $\disc_{t4}$, including
$5600$ acute triangles, $\disc_{t5}$, including
$4\times 5600 = 22400$ acute  triangles, 
$\disc_{t6}$, including
$16\times 5600 = 89600$ acute  triangles, $\disc_{q4}$, including
$1600$ rectangles (in fact, squares), $\disc_{q5}$, including
$4\times 1600 = 6400$ rectangles, $\disc_{q6}$, including
$25\times 1600 = 40000$ rectangles. 
For the triangular grids $\disc_{t4}$, $\disc_{t5}$ and $\disc_{t6}$,  
the points $\x_K$ have been located at the circumcenter of the triangles.
\begin{remark}\label{circum}
Choosing for  $\x_K$ the circumcenter of the triangle instead
of the center of gravity leads to an error about ten percent lower on 
the grids $\disc_{t4}$, $\disc_{t5}$ and $\disc_{t6}$.
\end{remark}
For the rectangular grids, the  points $\x_K$ have been located at  the center
of gravity of the control volumes.
We provide in the following table
the error $E_2 = \Vert u - \bar u\Vert_{L^2(\O)}$, as well as the minimum value
$u_{\min} = \min_{K\in\mesh} u_K$
and the maximum value $ u_{\max} =\max_{K\in\mesh} u_K$
of the approximate solution (note that the exact solution $\bar u$ 
varies between $0$ at the edges of $\O$ and $1$ at its center), using both schemes.

 \begin{center}
 \begin{tabular}{| c | c | c | c | c | c | c |}
\hline Grid
  & \begin{tabular}{c}MFE \\ $E_2$\end{tabular} 
  & \begin{tabular}{c}MFE \\ $u_{\min}$\end{tabular} 
  & \begin{tabular}{c}MFE \\ $u_{\max}$\end{tabular} 
  & \begin{tabular}{c}MFV \\ $E_2$\end{tabular} 
  & \begin{tabular}{c}MFV \\ $u_{\min}$\end{tabular} 
  & \begin{tabular}{c}MFV \\ $u_{\max}$\end{tabular} 
  \\
\hline
$\disc_{t4}$     &  1.53 &-1.32& 6.35 &  1.20  &-2.46&4.68  \\
$\disc_{t5}$   &  0.397  &-0.344& 2.20  &  0.315   &-0.633&1.99   \\
$\disc_{t6}$   &  0.101  &-0.0867& 1.20 &  0.0807   &-0.163&1.25   \\
$\disc_{q4}$    &  0.795 &-1.03& 2.62   &  0.000912   &0.000566& 0.997  \\
$\disc_{q5}$     &  0.200  &-0.259&1.38   &  0.000162  &0.000141& 0.999  \\
$\disc_{q6}$     &  0.0320  &-0.0415&1.06  &  0.0000202  & 0.0000229& 1.00  \\
\hline
\end{tabular}
\end{center}
 
\begin{figure}
\begin{center}
\begin{tabular}{ccc}
\includegraphics[width=0.25 \linewidth]{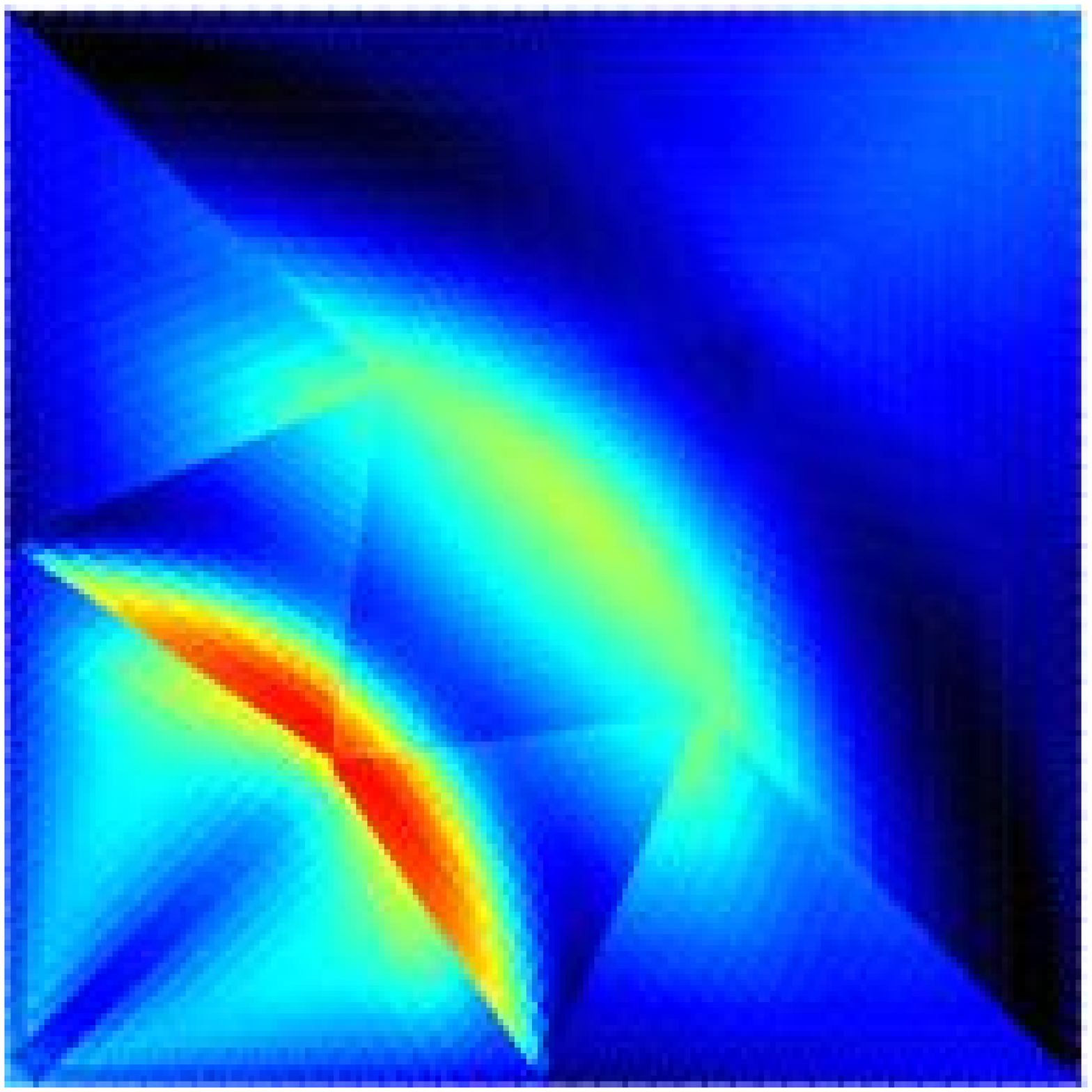}&
\includegraphics[width=0.25 \linewidth]{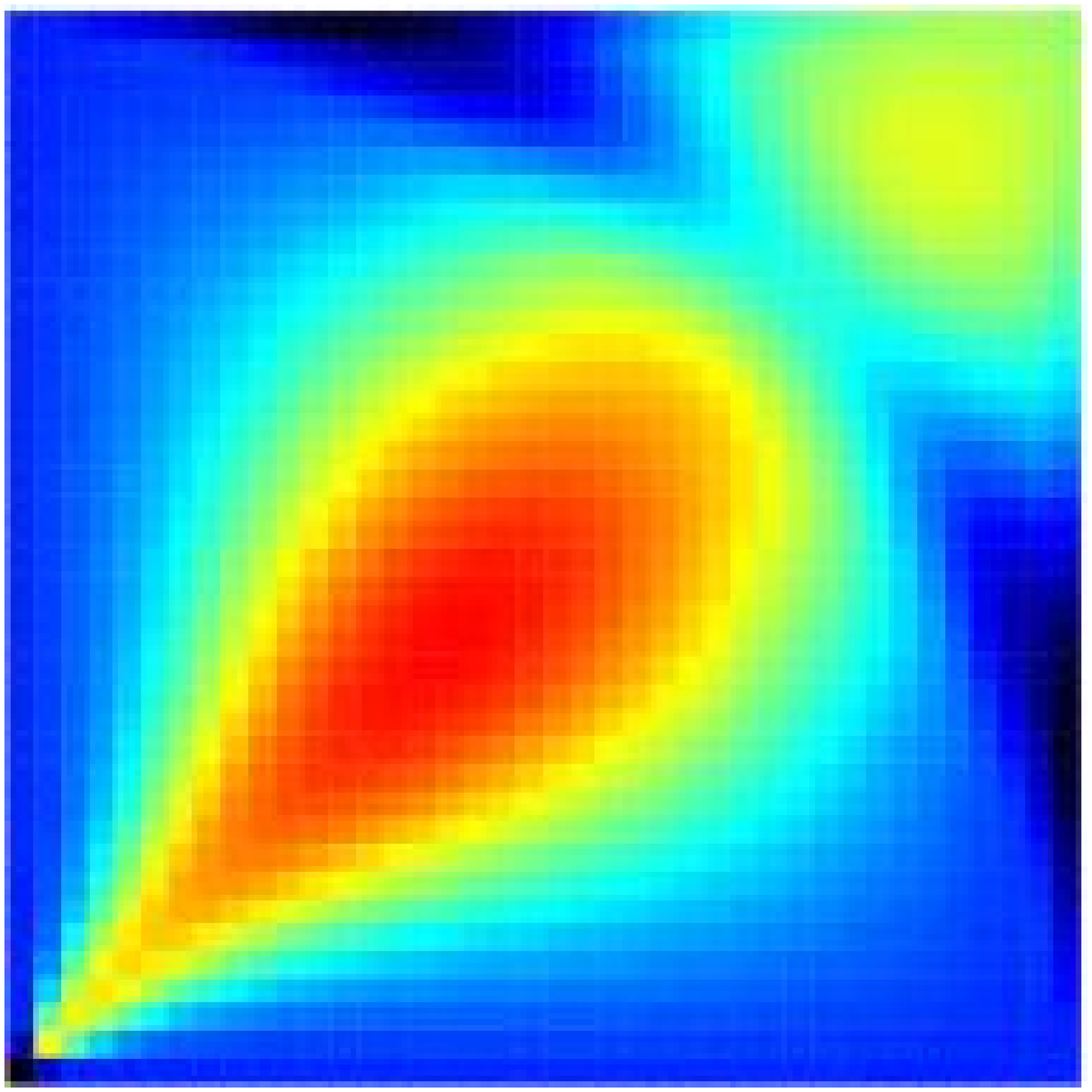}&
\includegraphics[width=0.25 \linewidth]{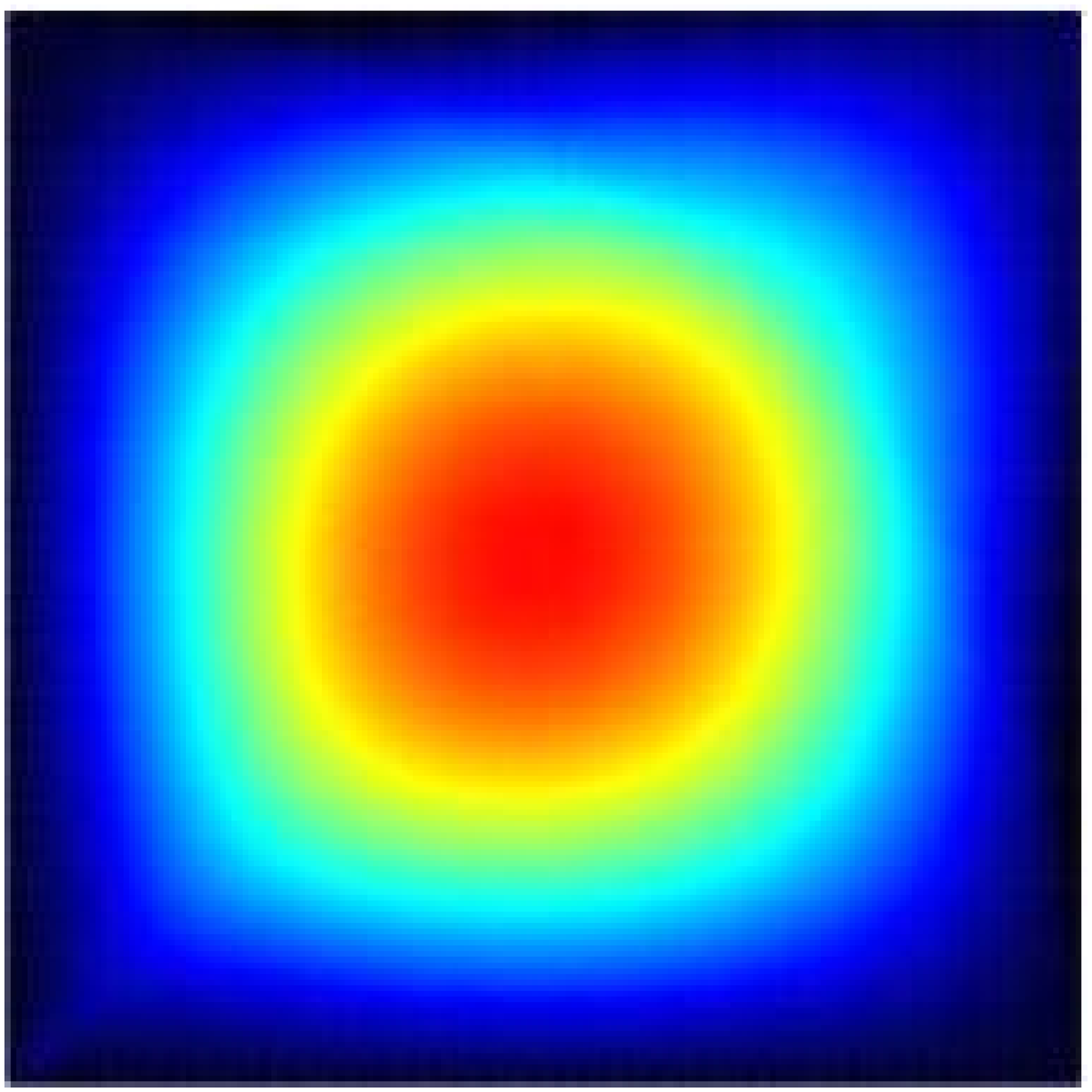}\\
MFE $\disc_{t4}$ & MFE $\disc_{q4}$ & MFE $\disc_{q6}$ \\ \\
\includegraphics[width=0.25 \linewidth]{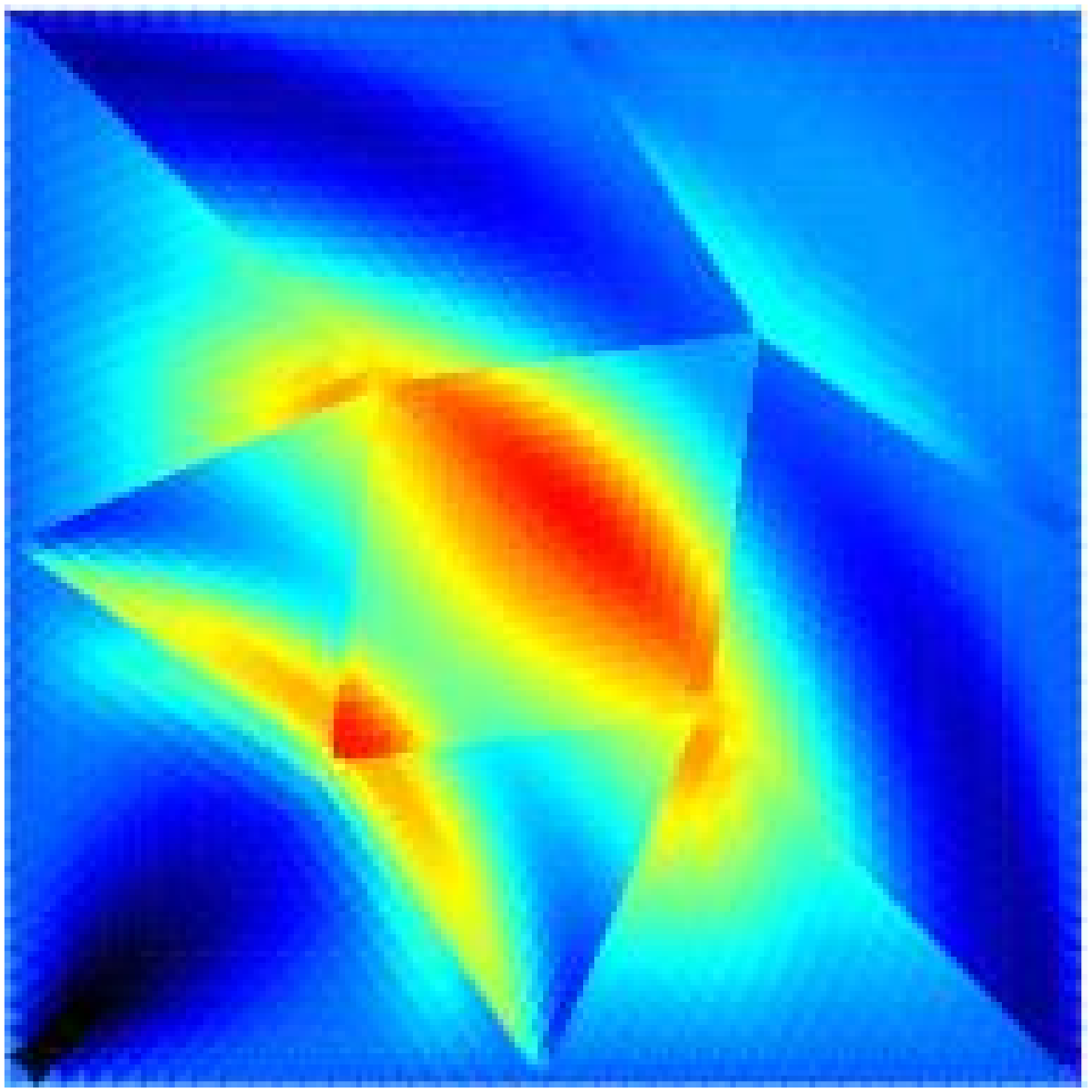}&
\includegraphics[width=0.25 \linewidth]{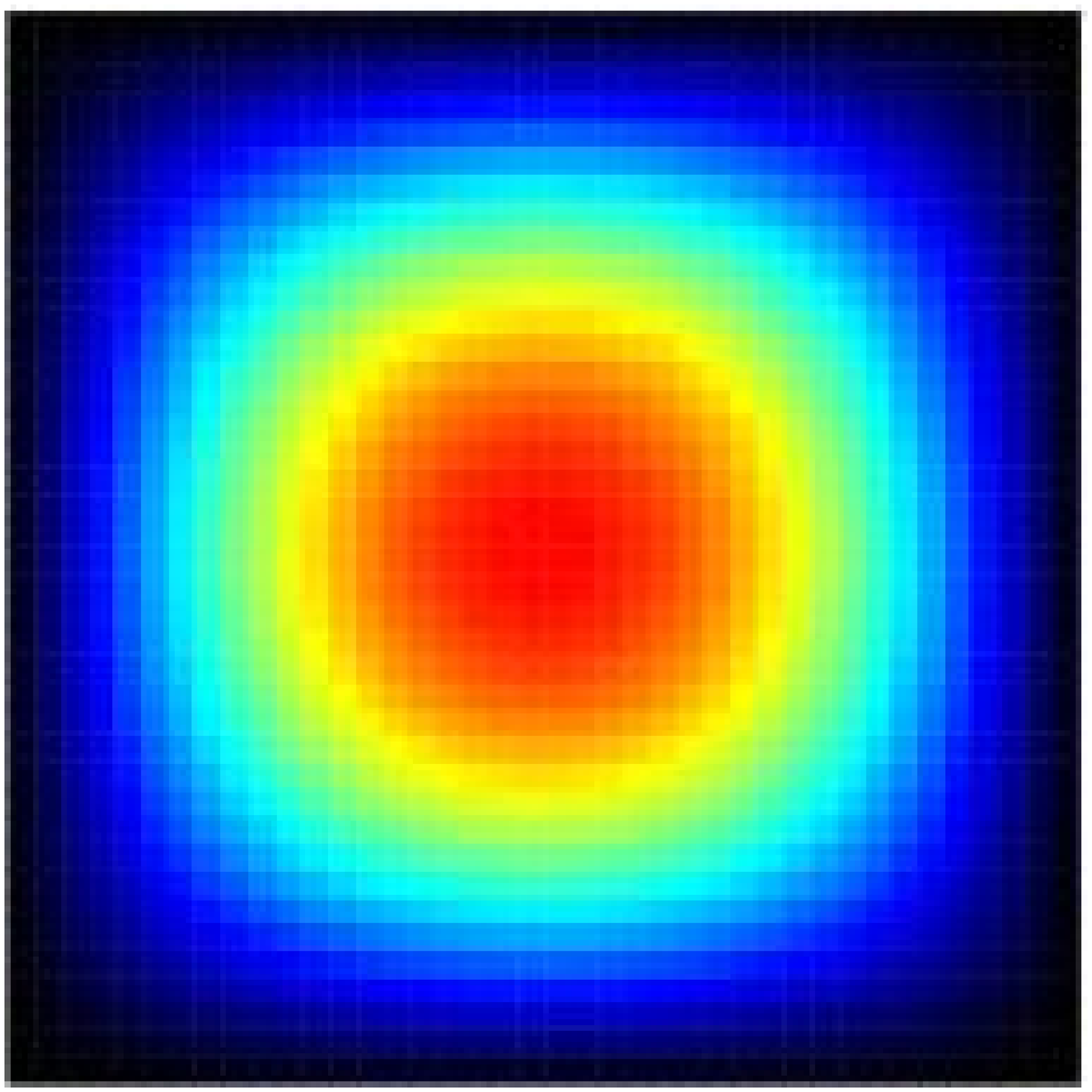}&
\includegraphics[width=0.25 \linewidth]{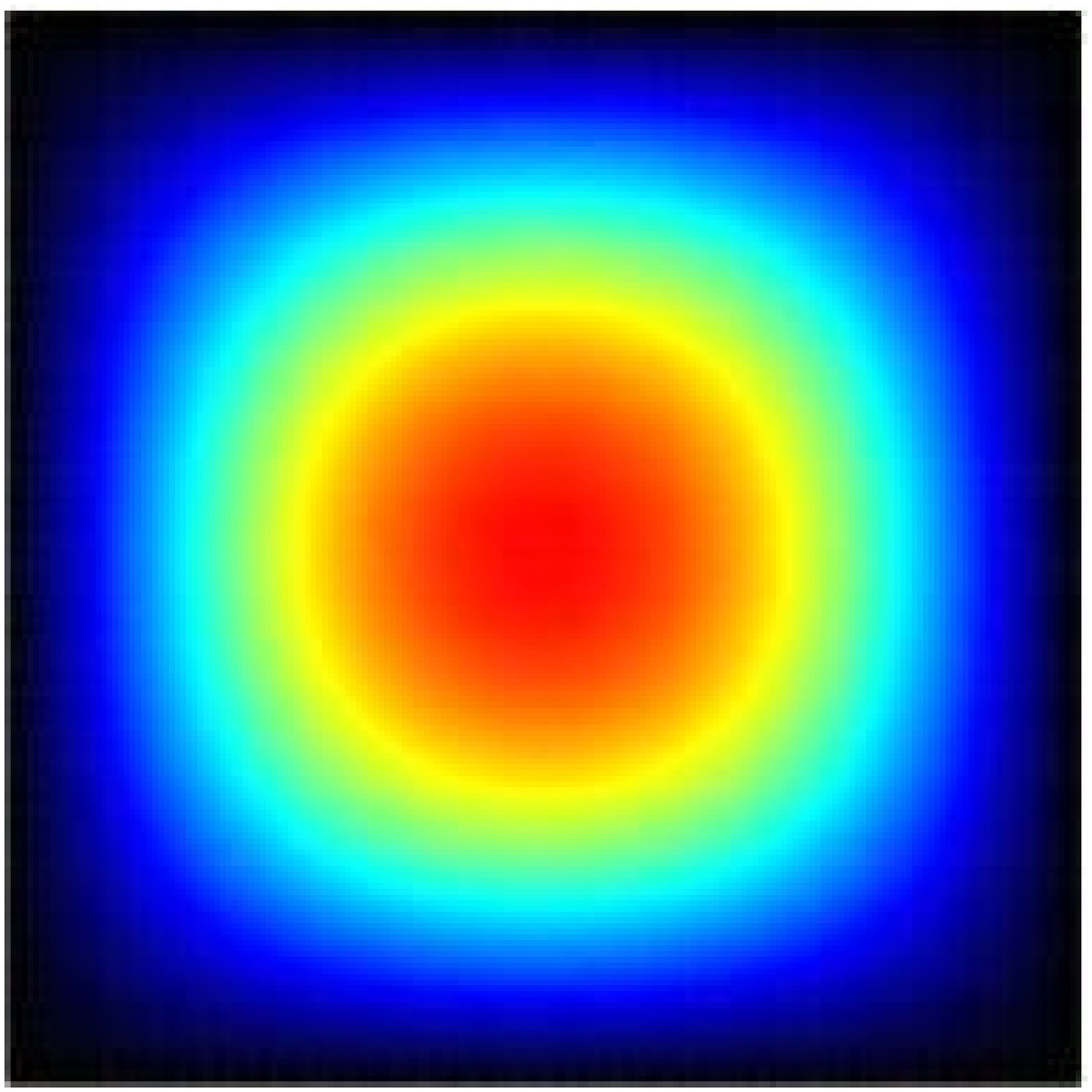}\\
MFV $\disc_{t4}$ & MFV $\disc_{q4}$ & MFV $\disc_{q6}$ 
\end{tabular}
\caption{Comparison of mixed finite element and mixed finite
volume solutions on Le Potier's test case.}
\label{cmplepo}
\end{center}
\end{figure}

These results show a surprisingly bad performance for the MFE and MFV schemes on triangular grids
(this was pointed out for the MFE scheme in \cite{lepot}). An order of convergence close to 2
is nevertheless observed for the $L^2(\O)$ norm of the error, with a very high multiplicative
constant. But this similarity between
both schemes does no longer hold on the other grids: 
on the regular rectangular grids (on which the MFE solution can be computed
using the classical RT basis), 
the MFV method provides accurate results
where the MFE scheme is far from the exact solution. 
Moreover,  in the case of the MFV scheme, the bounds on the approximate solution
are close to that of the exact solution.
These results are confirmed by Figure 
\ref{cmplepo}, where some of the numerical solutions considered in the above
table are plotted.

We now give in the following table the values
$E_2$,  $u_{\min}$ and $ u_{\max}$ in the case where the MFV scheme
is used on three irregular grids: the grid $\disc_{v}$ which is a Vorono\"{\i} tessellation
with 105 control volumes, the grid $\disc_{q1}$, already considered above, including
$16 + 144 + 49 + 25 = 234$ rectangles (in fact again, squares)
and the grid $\disc_{\sharp}$ with $400$ quadrangles, also considered above.

\begin{center}
\begin{tabular}{| c | c | c | c |}
\hline Grid
  & \begin{tabular}{c}MFV \\ $E_2$\end{tabular} 
  & \begin{tabular}{c}MFV \\ $u_{\min}$\end{tabular} 
  & \begin{tabular}{c}MFV \\ $u_{\max}$\end{tabular} 
  \\
\hline
$\disc_{v}$   &  0.0929 &0.0126& 0.980  \\
$\disc_{q1}$     &  0.0232 & 0.00259& 1.00  \\
$\disc_{\sharp}$    &  0.0217 & -0.00890& 0.999  \\
\hline
\end{tabular}
\end{center}

\begin{figure}
\begin{center}
\begin{tabular}{ccc}
\includegraphics[width=0.25 \linewidth]{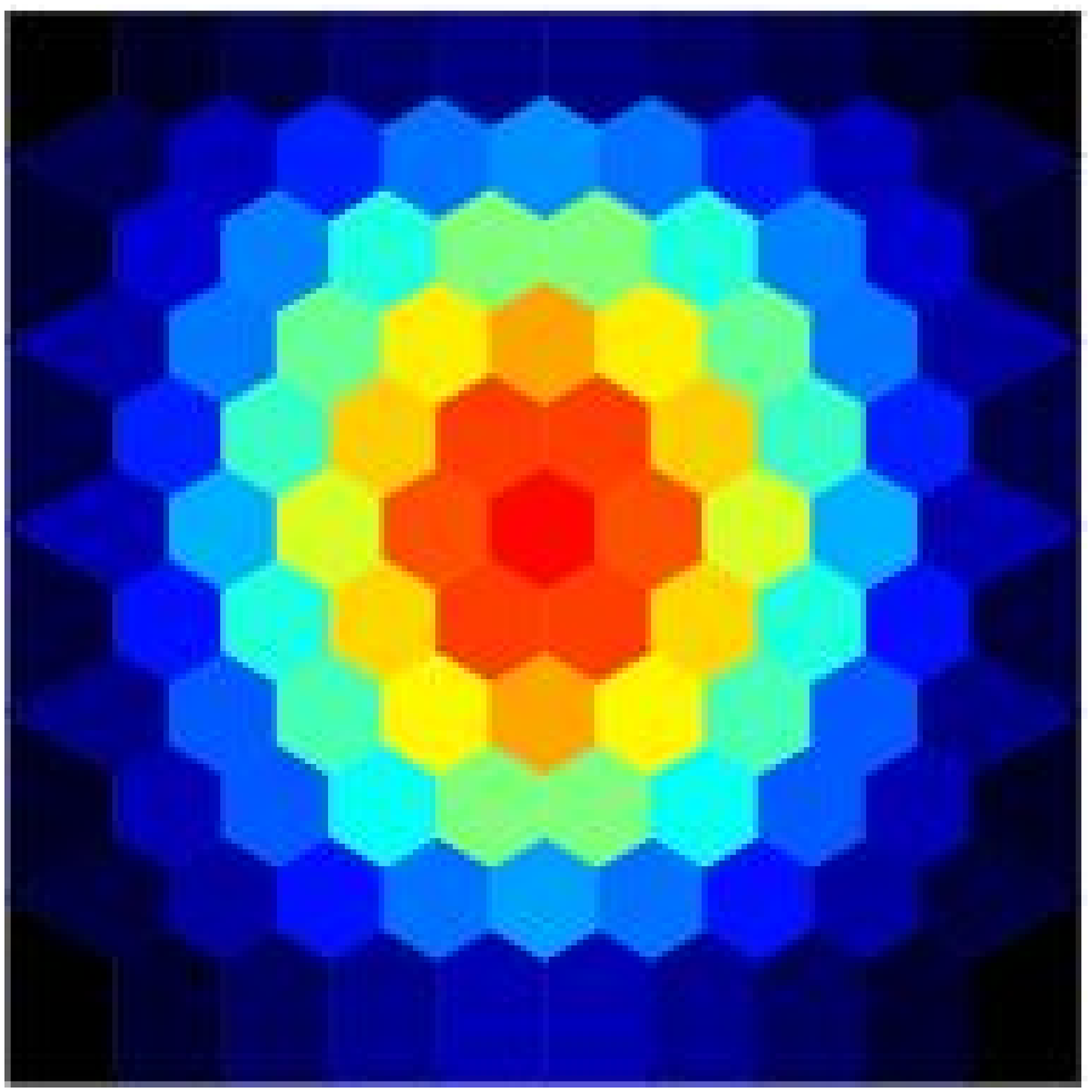}&
\includegraphics[width=0.25 \linewidth]{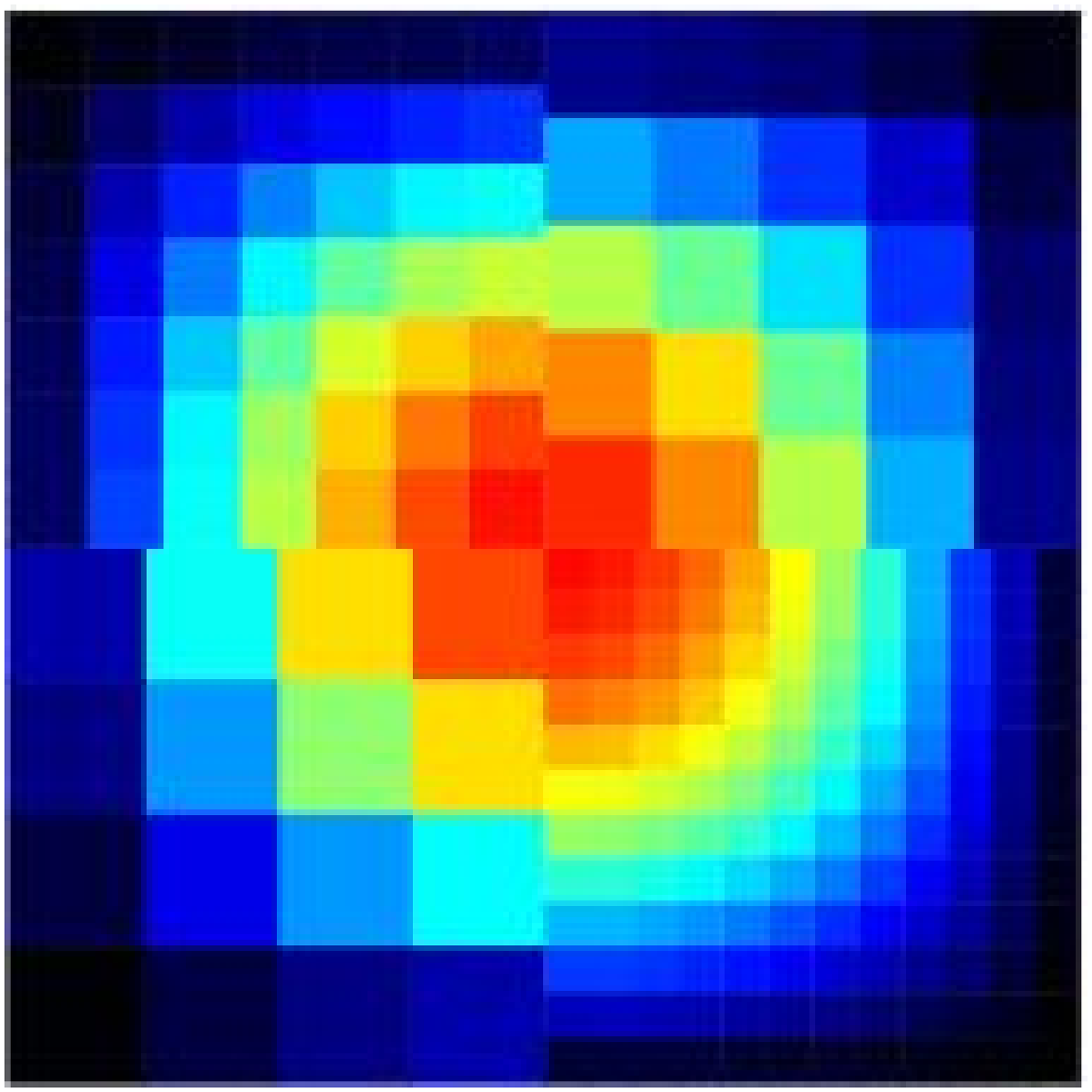}&
\includegraphics[width=0.25 \linewidth]{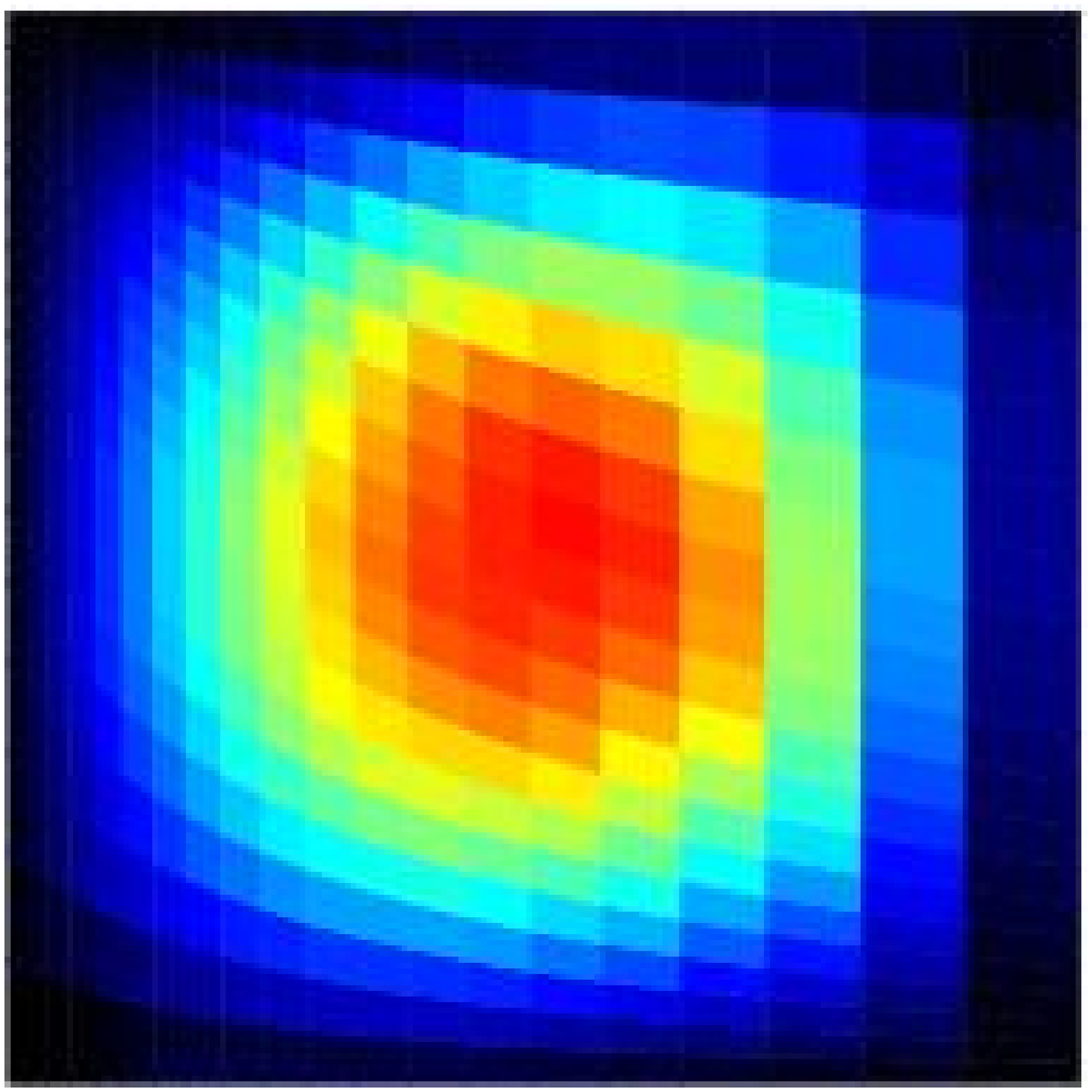} \\
MFV $\disc_v$ &MFV $\disc_{q1}$ & MFV $\disc_\sharp$
\end{tabular}
\caption{Solutions of Le Potier's test case using MFV on irregular grids.}
\label{griglepo}
\end{center}
\end{figure}

These results show an acceptable convergence, confirmed by Figure \ref{griglepo} in
which the corresponding approximate solutions are drawn.

\section{Appendix}

\subsection{Technical lemmas}

Lemma \ref{form-mag} justifies the link \refe{condconspen3} between the
approximate gradient and the approximate fluxes.

\begin{lemma}\label{form-mag}
Let $K$ be a non empty open convex polygonal set in $\R^d$.
For $\sigma\in\mathcal E_K$ (the edges of $K$, in the sense given in
Definition \ref{adisc}), we let $\mathbf{x}_\sigma$ be the
center of gravity of $\sigma$; we also denote $\n_{K,\sigma}$ the unit
normal to $\sigma$ outward to $K$. Then, for all vector $\mathbf{e}\in \R^d$
and for all point $\mathbf{x}_K\in \R^d$, we have
\[
\meas(K)\mathbf{e}=\sum_{\sigma\in\mathcal E_K} \meas(\sigma)\mathbf{e}\cdot\n_{K,\sigma}
(\mathbf{x}_\sigma-\mathbf{x}_K).
\]
\end{lemma}

\begin{proof}{.}

We denote by a superscript $i$ the $i$-th coordinate of vectors and
points in $\R^d$. By Stokes formula, we have
\[
\mcv\mathbf{e}^i=\int_K \div((x^i-\mathbf{x}_K^i)\mathbf{e})\dx
=\sum_{\sigma\in\mathcal E_K} \int_\sigma (x^i-\mathbf{x}_K^i)\mathbf{e}\cdot\n_{K,\sigma}\dfrontiere(x)
\]
and the proof is concluded since, by definition of the center of gravity, 
$\int_{\sigma}(x^i-\mathbf{x}_K^i)\dfrontiere(x)=\int_\sigma x^i\,d\gamma(x)-\meas(\sigma)\mathbf{x}_K^i
=\meas(\sigma)\mathbf{x}_\sigma^i-\meas(\sigma)\mathbf{x}_K^i$. \end{proof}

\medskip

The following lemma is quite similar to \cite[Lemma 7.2]{intervf},
but since the proof Lemma \ref{est-intervf} uses this result with slightly more general hypotheses
than in \cite{intervf}, we include the full proof of Lemma \ref{cop-col} for sake of completeness.

\begin{lemma}\label{cop-col}
Let $K$ be a non empty open polygonal convex set in $\R^d$.
Let $E$ be an affine hyperplane of $\R^d$ and $\sigma$ be a
non empty open subset of $E$ contained in $\partial K\cap E$.
We assume that there exists $\alpha>0$ and $\mathbf{p}_K\in K$
such that $B(\mathbf{p}_K,\alpha\diam(K))\subset K$.
We denote $\di_{K,\sigma}$ the convex hull of $\sigma$
and $\mathbf{p}_K$. Then there exists $\ctel{ze}$ only depending on $d$ and $\alpha$
such that, for all $v\in H^1(K)$,
\[
\left(\frac{1}{\meas(\di_{K,\sigma})}\int_{\di_{K,\sigma}} v(x)\dx
-\frac{1}{\meas(\sigma)}\int_{\sigma} v(\xi)\dfrontiere(\xi)\right)^2\le 
\frac{\cter{ze}{\rm dist}(\mathbf{p}_K,E)^2}{\meas(\di_{K,\sigma})}
\int_{\di_{K,\sigma}}|\nabla v(x)|^2\dx.
\]
\end{lemma}

\begin{proof}{.}

The regular functions being dense in $H^1(K)$ (since $K$ is convex), it is sufficient to prove
the lemma for $v\in C^1(\R^d)$. By translation and rotation, we can assume that $E=\{0\}\times\R^{d-1}$,
$\sigma=\{0\}\times\widetilde{\sigma}$ with $\widetilde{\sigma}\subset\R^{d-1}$
and that $\mathbf{p}_K=(p_1,0)$ with $p_1={\rm dist}(\mathbf{p}_K,E)$.

Notice that, since $K$ is convex and $\partial K\cap E$
contains a non empty open subset of $E$, $K$ is on one side of $E$.
In particular, $B(\mathbf{p}_K,\alpha\diam(K))$ is also on
one side of $E$ (it is contained in $K$) and
\begin{equation}
p_1={\rm dist}(\mathbf{p}_K,E)\ge \alpha\diam(K).
\label{hehe}\end{equation}

\medskip

For $a\in [0,p_1]$, we denote $\widetilde{\sigma}_a=\{z\in \R^{d-1}\;|\;(a,z)\in \di_{K,\sigma}\}$.
By definition, $(a,z)\in \di_{K,\sigma}$ if and only if there exists $t\in [0,1]$
and $y\in\widetilde{\sigma}$ such that $t(p_1,0)+(1-t)(0,y)=(a,z)$;
this is equivalent to $t=\frac{a}{p_1}$ and $z=(1-t)y=\left(1-\frac{a}{p_1}\right)y$.
Thus, $\widetilde{\sigma}_a=\left(1-\frac{a}{p_1}\right)\widetilde{\sigma}$.

\medskip

For all $y\in \widetilde{\sigma}$ and all $a\in [0,p_1]$, we have
\[
v(0,y)-v\left(a,\left(1-\frac{a}{p_1}\right)y\right)=
\int_0^1 \nabla v\left(ta,\left(1-t\frac{a}{p_1}\right)y\right)\cdot\left(-a,\frac{a}{p_1}y\right)\,{\rm d}t.
\]
Integrating on $y\in\widetilde{\sigma}$ and using the change
of variable $z=\left(1-\frac{a}{p_1}\right)y$, we find
\[
\int_{\sigma}v(\xi)\dfrontiere(\xi)-\frac{1}{\left(1-\frac{a}{p_1}\right)^{d-1}}\int_{\widetilde{\sigma}_a}v(a,z)\,{\rm d}z=
\int_{\widetilde{\sigma}}\int_0^1 \nabla v\left(ta,\left(1-t\frac{a}{p_1}\right)y\right)\cdot\left(-a,\frac{a}{p_1}y\right)
\,{\rm d}t{\rm d}y.
\]
Multiplying by $\left(1-\frac{a}{p_1}\right)^{d-1}$ and integrating on $a\in [0,p_1]$,
we obtain
\begin{eqnarray*}
\lefteqn{\int_{\sigma} v(\xi)\dfrontiere(\xi)\int_0^{p_1} \left(1-\frac{a}{p_1}\right)^{d-1}\,{\rm d}a
-\int_0^{p_1}\int_{\widetilde{\sigma}_a}v(a,z)\,{\rm d}z{\rm d}a}\\
&=&\int_0^{p_1}\left(1-\frac{a}{p_1}\right)^{d-1}\int_{\widetilde{\sigma}}\int_0^1 
\nabla v\left(ta,\left(1-t\frac{a}{p_1}\right)y\right)\cdot\left(-a,\frac{a}{p_1}y\right)\,{\rm d}t{\rm d}y{\rm d}a.
\end{eqnarray*}
But $\int_0^{p_1} \left(1-\frac{a}{p_1}\right)^{d-1}\,{\rm d}a=\frac{p_1}{d}$
and $\meas(\di_{K,\sigma})=\frac{\meas(\sigma)p_1}{d}$; therefore,
dividing by $\meas(\di_{K,\sigma})$, we find
\begin{eqnarray}
\lefteqn{\frac{1}{\meas(\sigma)}\int_{\sigma}v(\xi)\dfrontiere(\xi)
-\frac{1}{\meas(\di_{K,\sigma})}\int_{\di_{K,\sigma}}v(x)\,{\rm d}x}&&\nonumber\\
&=&\frac{1}{\meas(\di_{K,\sigma})}\int_0^{p_1}\left(1-\frac{a}{p_1}\right)^{d-1}\int_{\widetilde{\sigma}}\int_0^1 
\nabla v\left(ta,\left(1-t\frac{a}{p_1}\right)y\right)\cdot\left(-a,\frac{a}{p_1}y\right)\,{\rm d}t{\rm d}y{\rm d}a.
\label{intbord2}\end{eqnarray}

\medskip

For all $y\in\widetilde{\sigma}$, we have $|y|=|(0,y)|\le |(0,y)-\mathbf{p}_K|+|\mathbf{p}_K|\le
\diam(K)+p_1$ (because $(0,y)$ and $\mathbf{p_K}$ belong to $\overline{K}$).
By \refe{hehe}, this implies $|y|\le (\frac{1}{\alpha}+1)p_1$
and thus
\begin{eqnarray}
\lefteqn{\left|\int_0^{p_1}\left(1-\frac{a}{p_1}\right)^{d-1}\int_{\widetilde{\sigma}}\int_0^1 
\nabla v\left(ta,\left(1-t\frac{a}{p_1}\right)y\right)\cdot\left(-a,\frac{a}{p_1}y\right)\,{\rm d}t{\rm d}y{\rm d}a\right|}&&\nonumber\\
&\le& \ctel{dep}\int_0^{p_1}\left(1-\frac{a}{p_1}\right)^{d-1}\int_{\widetilde{\sigma}}\int_0^1 
\left|\nabla v\left(ta,\left(1-t\frac{a}{p_1}\right)y\right)\right|a\,{\rm d}t{\rm d}y{\rm d}a\nonumber\\
&\le& \cter{dep}\int_0^{p_1}\int_{\widetilde{\sigma}}\int_0^1 
\left|\nabla v\left(ta,\left(1-t\frac{a}{p_1}\right)y\right)\right|a\left(1-\frac{ta}{p_1}\right)^{d-1}\,{\rm d}t{\rm d}y{\rm d}a
\label{intbord3}\end{eqnarray}
where $\cter{dep}$ only depends on $\alpha$
(we have used the obvious fact that, for $t\in ]0,1[$, $1-\frac{a}{p_1}\le 1-\frac{ta}{p_1}$).
But, for all $a\in ]0,p_1[$, the change of variable 
\[
\varphi_a:(t,y)\in ]0,1[\times\widetilde{\sigma}\to
z=\left(ta,\left(1-t\frac{a}{p_1}\right)y\right)\in \varphi_a(]0,1[\times\widetilde{\sigma})
\]
has Jacobian determinant equal to $a\left(1-\frac{ta}{p_1}\right)^{d-1}$ and therefore
\[
\int_{\widetilde{\sigma}}\int_0^1 \left|\nabla v\left(ta,\left(1-t\frac{a}{p_1}\right)y\right)\right|
a\left(1-\frac{ta}{p_1}\right)^{d-1}\,{\rm d}t{\rm d}y=\int_{\varphi_a(]0,1[\times\widetilde{\sigma})}
|\nabla v(z)|\,{\rm d}z.
\]
Moreover, $(ta,(1-t\frac{a}{p_1})y)=\frac{ta}{p_1}(p_1,0)+(1-\frac{ta}{p_1})(0,y)$
with $\frac{ta}{p_1}\in ]0,1[$; hence, $\varphi_a(]0,1[\times\widetilde{\sigma})\subset \di_{K,\sigma}$
and we obtain
\[
\int_0^{p_1}\int_{\widetilde{\sigma}}\int_0^1 
\left|\nabla v\left(ta,\left(1-t\frac{a}{p_1}\right)y\right)\right|a\left(1-\frac{ta}{p_1}\right)^{d-1}\,{\rm d}t{\rm d}y{\rm d}a
\le p_1\int_{\di_{K,\sigma}} |\nabla v(z)|\,{\rm d}z.
\]
We introduce this inequality in \refe{intbord3} and use the resulting
estimate in \refe{intbord2} to obtain
\[
\left|\frac{1}{\meas(\di_{K,\sigma})}\int_{\di_{K,\sigma}} v(x)\,{\rm d}x
-\frac{1}{\meas(\sigma)}\int_{\sigma} v(\xi)\dfrontiere(\xi)\right|\le \frac{\cter{dep} p_1}{\meas(\di_{K,\sigma})}
\int_{\di_{K,\sigma}}|\nabla v(x)|\,{\rm d}x
\]
and the conclusion follows from the Cauchy-Schwarz inequality, recalling
that $p_1={\rm dist}(\mathbf{p}_K,E)$. \end{proof}

\medskip

\begin{lemma}\label{est-intervf}
Let $K$ be a non empty open polygonal convex set in $\R^d$ such that, for some $\alpha>0$, there exists
a ball of radius $\alpha \diam(K)$ contained in $K$. Let $E$ be an
affine hyperplane of $\R^d$ and $\sigma$ be a non empty open subset of $E$
contained in $\partial K\cap E$. Then there exists $\ctel{cst-je}$ only depending on $d$ and $\alpha$
such that, for all $v\in H^1(K)$,
\[
\left(\frac{1}{\mcv}\int_K v(x)\dx -\frac{1}{\meas(\sigma)}\int_\sigma v(x)\dfrontiere(x)\right)^2
\le \frac{\cter{cst-je}\diam(K)}{\meas(\sigma)}\int_K |\nabla v(x)|^2\dx.
\]
\end{lemma}

\begin{proof}{.}

Let $B(\mathbf{p}_K,\alpha\diam(K))\subset K$ and $\di_{K,\sigma}$ be the
convex hull of $\mathbf{p}_K$ and $\sigma$. By Lemma \ref{cop-col},
we have
\[
\left(\frac{1}{\meas(\di_{K,\sigma})}\int_{\di_{K,\sigma}} v(x)\dx -
\frac{1}{\meas(\sigma)}\int_\sigma v(x)\dfrontiere(x)\right)^2
\le \frac{\cter{ze}{\rm dist}(\mathbf{p}_K,E)^2}{\meas(\di_{K,\sigma})}\int_K |\nabla v(x)|^2\dx.
\]
But $\meas(\di_{K,\sigma})=\frac{\meas(\sigma){\rm dist}(\mathbf{p}_K,E)}{d}$ and
${\rm dist}(\mathbf{p}_K,E)\le {\rm dist}(\mathbf{p}_K,\sigma)\le \diam(K)$.
Therefore,
\begin{equation}
\left(\frac{1}{\meas(\di_{K,\sigma})}\int_{\di_{K,\sigma}} v(x)\dx -
\frac{1}{\meas(\sigma)}\int_\sigma v(x)\dfrontiere(x)\right)^2
\le \frac{\cter{ze}d\,\diam(K)}{\meas(\sigma)}\int_K |\nabla v(x)|^2\dx.
\label{un}\end{equation}

Using Lemma 7.1 in \cite{intervf}, we get $\ctel{cst-2}$ only depending on $d$ such that
\[
\left(\frac{1}{\meas(\di_{K,\sigma})}\int_{\di_{K,\sigma}}v(x)\dx
-\frac{1}{\mcv}\int_K v(x)\dx\right)^2
\le \frac{\cter{cst-2}\diam(K)^{d+2}}{\meas(\di_{K,\sigma})\mcv}\int_K|\nabla v(x)|^2\dx,
\]
which implies
\[
\left(\frac{1}{\meas(\di_{K,\sigma})}\int_{\di_{K,\sigma}}v(x)\dx
-\frac{1}{\mcv}\int_K v(x)\dx\right)^2
\le \frac{\cter{cst-2}d\,\diam(K)^{d+2}}{\meas(\sigma){\rm dist}(\mathbf{p}_K,E)\mcv}\int_K|\nabla v(x)|^2\dx.
\]
But, as in the proof of Lemma \ref{cop-col}, we have
${\rm dist}(\mathbf{p}_K,E)\ge \alpha\diam(K)$ (see \refe{hehe}).
Since $\mcv\ge \omega_d \alpha^d \diam(K)^d$, we deduce that
\begin{equation}
\left(\frac{1}{\meas(\di_{K,\sigma})}\int_{{\!\bigtriangleup\!}_{K,\sigma}}v(x)\dx
-\frac{1}{\mcv}\int_K v(x)\dx\right)^2
\le \frac{\cter{cst-2}d\,\diam(K)}{\omega_d\alpha^{d+1}\meas(\sigma)}\int_K|\nabla v(x)|^2\dx.
\label{deux}\end{equation}
The lemma follows from \refe{un} and \refe{deux}. \end{proof}

\subsection{Simplicial meshes}\label{sixdeux}

For some meshes, it is possible to completely drop the penalization on the
fluxes, that is to say to take $\nu_K=0$ in \refe{condconspen2}. This is for
example the case if each control volume $K$ of the mesh
is a simplex, i.e. if $K$ is the interior of the convex hull of 
$d+1$ points of $\R^d$ such that no affine hyperplane of $\R^d$ contains all of them
and if $\mbox{Card}(\edgescv) = d+1$.
In this situation, the following lemma is the key ingredient to the
study of the mixed finite volume scheme with $\nu_K=0$.

\begin{lemma}\label{estflu} 
Let us assume Assumptions \refe{hypomega}-\refe{hypfg}. Let $\disc$
be an admissible discretization of $\Omega$ in the sense of Definition \ref{adisc},
such that $\regul(\disc)\le \theta$ for some $\theta>0$ and $\mesh$ is made of simplicial
control volumes. Let $\mathbf{v}\in H_\disc^d$ and a family of real numbers 
$F=(F_{K,\sigma})_{K\in\mesh\,,\;\sigma\in\mathcal E_K}$ be given 
such that \refe{condconspen3} and \refe{mul3pen} hold. 
Then there exists  $\ctel{jjflu}$ only depending
on $d$, $\Omega$, $\Lambda$ and $\theta$ such that
\be
\sum_{K\in\mesh}\sum_{\sigma\in\mathcal E_K}\diam(K)^{2-d} F_{K,\sigma}^2
\le \cter{jjflu} (||f||_{L^2(\Omega)}^2 + ||\mathbf{v}||_{L^2(\Omega)^d}^2).
\label{eqestflu}\ee
\end{lemma}

\begin{proof}{.}

For $K\in\mesh$, let $A_K$ be the $(d+1)\times (d+1)$ matrix whose columns are 
$(1,\mathbf{x}_\sigma-\mathbf{x}_K)^T_{\sigma\in \mathcal E_K}$ (since $K$ is simplicial, 
it has $d+1$ edges and $A_K$
is indeed a square matrix). The equations \refe{condconspen3}-\refe{mul3pen} can be written $A_KF_K=E_K$,
where $F_K=(F_{K,\sigma})_{\sigma\in\mathcal E_K}$ and 
\[
E_K=\left(\begin{array}{c} -\int_K f(x)\dx \\ \mcv\Lambda_K\mathbf{v}_K\end{array}\right).
\]
We now want to estimate $||A_K^{-1}||$ and, in order to achieve this, we divide the rest 
of the proof in several steps.

\medskip

\emph{Step 1}: this step is devoted to allow the assumption $\diam(K)=1$ in Steps 2 and 3.

Let $K_0=\diam(K)^{-1}K$. Then $\mathbf{x}_{K,0}=\diam(K)^{-1}\mathbf{x}_K\in K_0$
and the barycenters of the edges of $K_0$ are $\mathbf{x}_{\sigma,0}=\diam(K)^{-1}\mathbf{x}_\sigma$.
Notice also that, if $\rho_{K,0}$ is the supremum of the radius of the balls
included in $K_0$, then 
\begin{equation}
\frac{1}{\rho_{K,0}}=\frac{\diam(K_0)}{\rho_{K,0}}=\frac{\diam(K)}{\rho_{K}}\le \regul(\disc)^{1/d}
\le \theta^{1/d}.
\label{ball}\end{equation}

Let $A_{K,0}$ be the $(d+1)\times(d+1)$ matrix corresponding to $K_0$, that is to say whose columns
are $(1,\mathbf{x}_{\sigma,0}-\mathbf{x}_{K,0})^T_{\sigma\in\mathcal E_K}=
(1,\diam(K)^{-1}(\mathbf{x}_\sigma-\mathbf{x}_K))^T_{\sigma\in \mathcal E_K}$. Since 
\[
A_K=\left( \begin{array}{cccc} 1 & 0 & \cdots & 0 \\ 0 & \diam(K) & \ddots & \vdots  \\ \vdots & \ddots &
\ddots & 0\\ 0&\cdots &0 & \diam(K)
\end{array} \right)A_{K,0},
\]
we have $||A_K^{-1}||\le \sup(1,\diam(K)^{-1})||A_{K,0}^{-1}||$. Hence, an estimate
on $||A_{K,0}^{-1}||$ gives an estimate on $||A_K^{-1}||$.

\medskip

\emph{Step 2}: estimate on $A_{K,0}$.

By \refe{ball}, $K_0$ contains a closed ball of radius $\frac 1 2 \theta^{-1/d}$. Up to a translation
(which does not change the vectors $\mathbf{x}_{\edge,0}-\mathbf{x}_{K,0}$, and hence
does not change $A_{K,0}$), we can assume that this ball is centered at $0$.
Since $\diam(K_0)=1$, we have then ${\overline B}(0,\frac 1 2 \theta^{-1/d})\subset K_0\subset {\overline B}(0,1)$.

Let $Z_\theta$ be the set of couples $(L,\mathbf{x}_L)$, where $L$ is a simplex
such that ${\overline B}(0,\frac 1 2 \theta^{-1/d})\subset \overline{L}\subset {\overline B}(0,1)$ and $x_L\in \overline{L}$.
Each simplex is defined by $d+1$ vertices in $\R^d$ so $Z_\theta$ can be
considered as a subset of $P=(\R^d)^{d+1}/S_{d+1}\times \R^d$,
where $S_{d+1}$ is the symmetric group acting on $(\R^d)^{d+1}$ by
permuting the vectors. As such, $Z_\theta$ is compact in $P$: it is
straightforward if we express the condition ``the adherence of a simplex contains
${\overline B}(0,\frac 1 2 \theta^{-1/d})$'' as ``any point of ${\overline B}(0,\frac 1 2 \theta^{-1/d})$
is a convex combination of the vertices of the simplex'', which is a closed condition with respect to
the vertices of the simplex.

For $(L,\mathbf{x}_L)\in Z_\theta$, let $M(L,\mathbf{x}_L)$ be the set of
$(d+1)\times(d+1)$ matrices whose columns are, up to permutations, 
$(1,\mathbf{x}_\sigma-\mathbf{x}_L)^T_{\sigma\in\mathcal E_L}$ ($\mathcal E_L$
being the set of edges of $L$ and $\mathbf{x}_\sigma$ being the barycenter of $\sigma$).
$M(L,\mathbf{x}_L)$ can be considered as an element of $M_{d+1}(\R)/S_{d+1}$ ($S_{d+1}$ acting by permuting
the columns) and the application $(L,\mathbf{x}_L)\in Z_\theta \to M(L,\mathbf{x}_L)\in M_{d+1}(\R)/S_{d+1}$
is continuous: to see this, just recall that the barycenter
of an edge $\sigma\in\mathcal E_L$ is $\mathbf{x}_{\sigma}=\frac{1}{d}\sum_{i=1}^d \mathbf{x}_i$,
where $\mathbf{x}_i$ are the vertices of $\sigma$ (i.e. all vertices but one of
$L$).

If $(L,\mathbf{x}_L)\in Z_\theta$, all the matrices of $M(L,\mathbf{x}_L)$ are invertible.
Indeed, assume that such a matrix has a non-trivial element $(\lambda_1,\ldots,\lambda_{d+1})$
in its kernel; this leads (denoting $(\sigma_1,\ldots,\sigma_{d+1})$
the edges of $L$) to $\sum_{i=1}^{d+1} \lambda_i=0$ and $\sum_{i=1}^{d+1} \lambda_i(\mathbf{x}_{\sigma_i}
-\mathbf{x}_L)=\sum_{i=1}^{d+1} \lambda_i\mathbf{x}_{\sigma_i}=0$.
Assuming $\lambda_{d+1}\not= 0$, we then can write 
$\mathbf{x}_{\sigma_{d+1}}=\sum_{i=1}^d \mu_i \mathbf{x}_{\sigma_i}$ with
$\sum_{i=1}^d \mu_i=1$ (since $\mu_i=-\lambda_i/\lambda_{d+1}$). This means that
$\mathbf{x}_{\sigma_{d+1}}$ is in the affine hyperplane $\mathcal H$ generated by the
other barycenters of edges. Note that $\mathcal H$ is parallel to $\sigma_{d+1}$
(this is a straightforward consequence of 
Thales' theorem at the vertex which does not belong to $\sigma_{d+1}$,
and of the fact that the barycenters $(\mathbf{x}_{\sigma_1},\ldots,\mathbf{x}_{\sigma_d})$
of the edges are in fact the barycenters of the vertices of the
corresponding edge). Therefore  $\mathcal H$  contains the whole edge  $\sigma_{d+1}$, because it
contains  $\mathbf{x}_{\sigma_{d+1}}\in \sigma_{d+1}$. Let $\mathbf{a}$
be the vertex of $L$ which does not belong to $\sigma_{d+1}$;
$\mathbf{a}$ belongs to $\sigma_1$ and we denote $(\mathbf{b}_1,\ldots,\mathbf{b}_{d-1})$
the other vertices of $\sigma_1$ (which also belong to $\sigma_{d+1}$). We have $\mathbf{x}_{\sigma_1}=\frac{1}{d}(\mathbf{a}
+\sum_{i=1}^{d-1}\mathbf{b}_i)$, and therefore $\mathbf{a}=d\mathbf{x}_{\sigma_1}-\sum_{i=1}^{d-1}
\mathbf{b}_i$; but $d-\sum_{i=1}^{d-1} 1 = 1$ and thus $\mathbf{a}$
belongs to the affine hyperplane generated by $(\mathbf{x}_{\sigma_1},\mathbf{b}_1,\ldots,
\mathbf{b}_{d-1})$. Since all these points belong to $\mathcal H$, we
have $\mathbf{a}\in \mathcal H$ and, since $\sigma_{d+1}\subset \mathcal H$,
all the vertices of $L$ in fact belong to $\mathcal H$; $L$ is thus contained in an hyperplane,
which is a contradiction with the fact that it contains a non-trivial ball.
Thus, for $(L,\mathbf{x}_L)\in Z_\theta$, $M(L,\mathbf{x}_L)$ is in fact an element
of $Gl_{d+1}(\R)/S_{d+1}$.

The inversion ${\rm inv}:Gl_{d+1}(\R)\to Gl_{d+1}(\R)$
is continuous; hence, $||{\rm inv}(\cdot)||:Gl_{d+1}(\R)\to \R$ is also continuous.
Permuting the columns of a matrix comes down to permuting the lines of its inverse,
which does not change the norm; therefore $||{\rm inv}(\cdot)||:Gl_{d+1}(\R)/S_{d+1}\to \R$
is well defined and also continuous.

We can now conclude this step. The application $Z_{\theta}\to Gl_{d+1}(\R)/S_{d+1}\to \R$ defined
by $(L,\mathbf{x}_L)\to M(L,\mathbf{x}_L)\to ||{\rm inv}(M(L,\mathbf{x}_L))||$ is continuous.
Since $Z_\theta$ is compact, this application is bounded by some $\ctel{pfff}$ only
depending on $d$ and $\theta$. As $(K_0,\mathbf{x}_{K,0})\in Z_\theta$,
this shows that $||A_{K,0}^{-1}||\le \cter{pfff}$.

\medskip

\emph{Step 3}: conclusion.

Using the preceding steps, we find $||F_K||\le ||A_K^{-1}||\,||E_K||
\le \cter{pfff}\sup(1,\diam(K)^{-1})||E_K||$. Hence,
\[
\sum_{K\in\mesh}\diam(K)^{2-d} ||F_K||^2 \le \cter{pfff}^2\sup(\diam(\O)^2,1)\sum_{K\in\mesh} \diam(K)^{-d}||E_K||^2.
\]
But $||E_K||^2\le \mcv \int_K |f(x)|^2\dx+\ctel{pf}\mcv^2 |\mathbf{v}_K|^2$
with $\cter{pf}$ only depending on $\Lambda$. 
Since $\mcv\le \omega_d \diam(K)^d$, this concludes the proof of 
\refe{eqestflu}. \end{proof}

\medskip

Let us now consider (\refe{condconspen2},\refe{condconspen1},\refe{condconspen3},\refe{mul3pen}) with $\nu_K=0$;
notice that the results of Section \ref{sec-discspace} still hold in this
situation. 

Equation \refe{convfort2} gives, if $\nu_K=0$, an estimate on $\mathbf{v}$ in
$L^2(\Omega)^d$ which, thanks to Lemma \ref{estflu}, translates into an estimate
on the fluxes (this estimate replaces the one obtained before thanks to the penalization),
provided that the control volumes are simplicial. This gives, as in Remark \ref{proofth}, existence and
uniqueness of a solution to the non-penalized mixed finite volume scheme (i.e. 
(\refe{condconspen2},\refe{condconspen1},\refe{condconspen3},\refe{mul3pen}) with $\nu_K=0$).
{}From the estimate on the fluxes, it is straightforward to see that the term $\termr{convpp}$ in the
proof of Theorem \ref{convpen} still tends to $0$ as $\size(\disc)\to 0$.
Hence, in the case of simplicial control volumes, the solution to
the mixed finite volume scheme
(\refe{condconspen2},\refe{condconspen1},\refe{condconspen3},\refe{mul3pen})
with $\nu_K=0$ still converges toward the weak solution of \refe{ellgen}.

It is also quite easy to establish, in this situation, error estimates in the case
of smooth data $\Lambda$ and $\bar u$; these estimates are in fact quite better than
the ones of Theorem \ref{esterrorpen}: we can prove that
\[
\Vert \mathbf{v} - \grad \bar u\Vert_{L^2(\O)^d} \le \ctel{cerr0}
\size(\disc)\qquad\mbox{and}\qquad
\Vert u - \bar u\Vert_{L^2(\O)} \le \cter{cerr0} \size(\disc).
\]
To obtain such rates of convergence, one must simply bound $\termr{penerr3}$
in \refe{dring} by using Lemma \ref{estflu} with $F=\widehat{F}$, $\mathbf{v}=\widehat{\mathbf{v}}$
and $f=0$.

\begin{remark}\label{expefm}
In the particular case where $\disc$ is made of simplicial
control volumes,  and, for all $K\in\mesh$, $\nu_K=0$ and $\x_K$ is the
center of gravity of $K$, then the solution  $(u,{\bf v},F)$ of 
(\refe{condconspen2},\refe{condconspen1},\refe{condconspen3},\refe{mul3pen})
is also the solution of the following generalization of the 
expanded mixed finite element scheme \cite{chen}:
find $(u,{\bf v},{\bf w}=\sum_{K\in\mesh}\sum_{\edge\in\edgescv} F_{K,\edge} {\bf W}_{K,\edge})
\in H_\disc\times H_\disc^d\times RT^0$  ($RT^0$ denotes here
the lowest degree Raviart-Thomas basis $({\bf W}_\edge)_{\edge\in\edges}$ on the mesh $\mesh$,
such that, choosing for an internal edge $\edge = K|L$ the orientation
from $K$ to $L$, then ${\bf W}_\edge$ restricted to $K$
is ${\bf W}_{K,\edge}$ and   ${\bf W}_\edge$ restricted to $L$ is $-{\bf W}_{L,\edge}$ --- note
that ${\bf w}\in RT^0$ thanks to \refe{condconspen1}) such that
\[
\int_\O \Lambda(x)\mathbf{v}(x)\cdot \mathbf{v}'(x) \dx = \int_\O
{\bf w}(x)\cdot \mathbf{v}'(x) \dx, \ \forall {\bf v}'\in H_\disc^d,\]
which gives \refe{condconspen3},
\[
\int_\O \mathbf{v}(x)\cdot {\bf w}'(x) \dx + \int_\O u(x)\div {\bf w}'(x) \dx = 0,
\forall {\bf w}'\in RT^0,
\]
which gives \refe{condconspen2} with $\nu_K=0$, and
\[
- \int_\O u'(x) \div {\bf w}(x) \dx = \int_\O u'(x) f(x) \dx,\ \forall u'\in H_\disc,
\]
which gives \refe{mul3pen}.
This formulation differs from that of  \cite{chen}, 
in which the restrictions of ${\bf v}$ and  ${\bf w}$ on each control volume must belong to the same space.
The proof of convergence of the mixed finite volume scheme therefore gives at the same time
that of this particular version of the expanded mixed finite element scheme.
\end{remark}

\end{document}